\documentclass[11pt, reqno]{amsart}

\usepackage[margin=1in]{geometry}
\usepackage{latexsym}
\usepackage{graphicx}
\usepackage{xypic}
\usepackage{hyperref}
\usepackage{amsmath,amssymb,amsthm,amsfonts,hyperref}
\usepackage{enumerate,enumitem}
\usepackage{subcaption}
\usepackage{tikz-cd}
\usepackage{lineno}
\usepackage{minibox}
\usepackage[ruled,linesnumbered]{algorithm2e}
\usepackage{appendix}
\allowdisplaybreaks
\hypersetup{%
    colorlinks,
    linkcolor={red!50!black},
    citecolor={green!50!black},
    urlcolor={blue!50!black}
}

\newtheorem{theorem}{Theorem}
\theoremstyle{plain}

\newtheorem{claim}[theorem]{Claim}
\newtheorem{question}{Question}

\newtheorem{construction}[theorem]{Construction}

\newtheorem{corollary}[theorem]{Corollary}

\newtheorem{definition}[theorem]{Definition}

\newtheorem{fact}[theorem]{Fact}
\newtheorem{lemma}[theorem]{Lemma}

\newtheorem{proposition}[theorem]{Proposition}
\newtheorem{prop}[theorem]{Proposition}

\numberwithin{equation}{section}
\numberwithin{theorem}{section}

\def \e{\varepsilon}

\def \a{\alpha}
\def\eps{\varepsilon}

\def \bfi{\mathbf{i}}

\def \bfv{\mathbf{v}}

\def\a{\alpha}

\def\d{\delta}
\def\g{\gamma}

\def\COMMENT#1{}
\let\COMMENT=\footnote
\newcommand{\red}{\textcolor{red}}

\begin{document}

\title{A parameterized algorithm for $K_r$-factors in graphs of high minimum degree}
\author{Luyining Gan}
\thanks{L. Gan was partially supported by the National Natural Science Foundation of China (12401446) and the Fundamental Research Funds for the Central Universities}
\address{School of Mathematical Sciences and Key Laboratory of Mathematics and Information Networks\\ Beijing University of Posts and Telecommunications\\ Beijing\\ China}
\email{elainegan@bupt.edu.cn}

\author{Jie Han}
\thanks{J. Han was partially supported by the National Natural Science Foundation of China (12371341)}
\address{School of Mathematics and Statistics and Center for Applied Mathematics\\ Beijing Institute of Technology\\ Beijing\\ China}
\email{han.jie@bit.edu.cn}

\author{Jie Hu}
\thanks{J. Hu was partially supported by the National Natural Science Foundation of China (12301443)}
\address{School of Mathematical Sciences and LPMC\\ Nankai University\\ Tianjin\\ China}
\email{hujie@nankai.edu.cn}

\begin{abstract}
A $K_r$-factor of a graph $G$ is a collection of vertex-disjoint $r$-cliques covering $V(G)$.
We prove the following algorithmic version of the classical Hajnal--Szemer\'edi Theorem in graph theory, when $r$ is considered as a constant.
Given $r, c, n\in \mathbb{N}$ such that $n\in r\mathbb N$, let $G$ be an $n$-vertex graph with minimum degree at least $(1-1/r)n - c$.
Then there is an algorithm with running time $2^{c^{O(1)}} n^{O(1)}$ that outputs either a $K_r$-factor of $G$ or a certificate showing that none exists, namely, this problem is fixed-parameter tractable in $c$.
On the other hand, it is known that if $c = n^{\varepsilon}$ for fixed $\varepsilon \in (0,1)$, the problem is \texttt{NP-C}. 
By taking the complement, our result yields a similar result on the equitable $\Delta$-colorings of graphs of maximum degree $\Delta+c$, for $\Delta\in [n/r, n/(r-1)]$.

We indeed establish characterization theorems for this problem, showing that the existence of a $K_r$-factor is equivalent to the existence of certain class of $K_r$-tilings of size $o(n)$, whose existence can be searched by the color-coding technique developed by Alon--Yuster--Zwick.

\end{abstract}
\date{\today}
\maketitle

\section{Introduction}

In this paper we concern about the decision/search problem for the existence of certain subgraphs.
We are interested in the target subgraphs whose order grows with the order of the host graph, as otherwise brute force search is enough.
The most prominent example is the perfect matching problem, famously shown by Edmonds Blossom Algorithm \cite{Edmonds} to be polynomial-time solvable.
However, Hell and Kirkpatrick~\cite{HeKi} showed that this is a rare example for spanning subgraphs -- they showed that the decision problem for essentially all other spanning subgraphs in general graphs is \texttt{NP-C}\footnote{To be more precise, they showed that the $F$-factor problem is \texttt{NP-C} if $F$ has a component on at least three vertices.}.
Therefore, it is natural to consider the restriction of the problem to certain class of host graphs so that the problem is tractable.

A classical line of research in extremal graph theory concerns about the smallest minimum degree condition or edge density condition that forces the existence of certain subgraphs.
This has been a rich and active research area since Dirac's theorem on Hamiltonian cycles~\cite{Dirac} proved in 1952.
The decision problems restricted to this class of graphs (with large minimum degree) are trivial but the search problems for general subgraphs are nice and challenging problems.

A natural next step to these results is to study the decision problem in the graphs where the studied parameter is \emph{slightly off the threshold}.
In particular, restricted to the minimum degree case, this is when the minimum degree is strictly below the threshold.
Indeed, this falls into the category of \emph{the above/below guarantee parameterization} introduced by Mahajan, Raman and Sikdar~\cite{MAHAJAN2009137}, which was successfully applied in the study of several fundamental problems in parameterized complexity and kernelization. 
For illustrative examples and a comprehensive review, see the recent excellent survey by Gutin and Mnich~\cite{GuMn}.

In particular, a classical theorem of Dirac~\cite{Dirac} back to 1952 says that every $n$-vertex 2-connected graph $G$ with minimum degree $\delta(G)\ge 2$ contains a cycle of length at least $\min\{2\delta(G), n\}$.
Fomin, Golovach, Sagunov, and Simonov~\cite{Fomin22} in 2022 proved an algorithmic extension of this result by showing that given $c\in \mathbb N$, there is an \texttt{FPT} algorithm parameterized in $c$ that determines whether such $G$ contains a cycle of length $\min\{2\delta(G)+c, n\}$, that is, the running time of the algorithm is $2^{O(c)}n^{O(1)}$ \footnote{In general, a decision problem with parameter $c$ is in \texttt{FPT} parameterized in $c$ if it is solvable in time $f(c)\cdot n^{O(1)}$ for some computable function $f$.}.
A similar result on the Dirac's theorem on Hamiltonicity is obtained earlier by Jansen, Kozma, and Nederlof~\cite{JKN} in 2019.
In view of these results, a general question can be stated as follows.

\begin{question}
\label{que}
Given $k>0$ and a graph theory theorem which says that every graph $G$ with $\delta(G) \ge \delta$ satisfies an (increasing) property $\mathcal P$.
How difficult is it to decide whether a graph $G$ with $\delta(G) \ge \delta - k$ has property $\mathcal P$, namely, $G\in \mathcal P$?
\end{question}

Although the minimum degree condition is probably the most studied parameter, it is clear that one could replace it in the question above by average degree, maximum degree or other graph parameters.

In this paper we take a further step towards Question~\ref{que} by studying an algorithmic extension of the celebrated Hajnal--Szemer\'edi Theorem.
To state the theorem we first need some definitions.
Given two graphs $F$ and $G$, an \emph{$F$-tiling} (also called $F$-packing, $F$-matching in other contexts) in $G$ is a collection of vertex-disjoint copies of $F$ in $G$.
The number of copies of $F$ is called the \emph{size} of the $F$-tiling.
An $F$-tiling is called \emph{perfect} if it covers all the vertices of $G$.
Perfect $F$-tilings are also referred to as \emph{$F$-factors}.

The following theorem of Hajnal and Szemer\'edi back to 1970 determines the sharp minimum degree condition to guarantee a $K_r$-factor in a graph.

\begin{theorem}~\cite[Hajnal--Szemer\'edi]{HSze} \label{thm:HSze}
Every graph $G$ whose order $n$ is divisible by $r$ and whose minimum degree satisfies $\d(G)\ge(1-1/r)n$ contains a $K_r$-factor. 
\end{theorem}

An \emph{equitable $k$-coloring} of a graph $G$ is a proper $k$-coloring, in which any two color classes differ in size by at most one. 
By considering the complement graph Theorem~\ref{thm:HSze} gives that every graph with maximum degree at most $\Delta$ has an equitable $(\Delta+1)$-coloring.
Indeed, Kierstead and Kostochka \cite{Kierstead-Kostochka-2008} gave a simpler proof of Theorem~\ref{thm:HSze} by proving the equitable coloring version.

Following Question~\ref{que}, the main goal of this paper is to study the following algorithmic problem stated in two equivalent forms.

\begin{center}
\minibox[frame]{
\textsc{Equitable $\Delta$-Coloring In Graphs Below Maximum Degree} \quad \hfill \textbf{Parameter}: $c$ \\
\textbf{Input:} $n$-vertex graph $G$ with maximum degree $\Delta + c$. \\
\textbf{Question:} Find an equitable $\Delta$-coloring in $G$ or report that it does not exist.
}
\end{center}

\begin{center}
\minibox[frame]{
$K_r$-\textsc{Factor In Graphs With High Minimum Degree} \quad \hfill \textbf{Parameter}: $c$ \\
\textbf{Input:} $n$-vertex graph $G$ with minimum degree $(1-1/r)n - c$. \\
\textbf{Question:} Find a $K_r$-factor in $G$ or report that it does not exist.
}
\end{center}


Our main result deals with this problem when $r$ is considered as a constant\footnote{Therefore, throughout this paper, when we use the big-$O$ notation, the implicit constant depends only on $r$, not $c$ or $n$.}.

\begin{theorem}
\label{thm:main}
Given $r, c, n\in \mathbb N$ with $n\in r\mathbb N$ and an $n$-vertex graph $G$ with $\d(G) \ge (1-\tfrac1r)n - c$, there is an algorithm with running time $n^{O(1)}2^{O(c^{{r}^{r+4}})}$ that decides whether $G$ has a $K_r$-factor, that is, the problem is \texttt{FPT} parameterized in $c$.
Moreover, the algorithm either returns a $K_r$-factor or a certificate showing that none exists.
In particular, for $r=3$, the algorithm runs in time $n^{O(1)}2^{O(c)}$.
\end{theorem}




For the sharpness of the degree condition in Theorem~\ref{thm:main}, when $c=0$, Theorem~\ref{thm:HSze} says that such $G$ must contain a $K_r$-factor and thus the decision problem is trivially in \texttt{P}.
Then Kierstead, Kostochka, Mydlarz and Szemer\'edi~\cite{KKMSz} gave a version of Theorem~\ref{thm:HSze} which yields a polynomial-time algorithm that produces a $K_r$-factor.
K\"uhn and Osthus~\cite{KuOs06soda} showed that if $c=\omega(n)$, then the decision problem of such $G$ containing a $K_r$-factor is \texttt{NP-C}, by reducing the problem to the $K_r$-factor problem in general graphs (with no minimum degree assumption).
We remark that this can be extended to that $c$ is any polynomial function of $n$ by the same proof.
\begin{proposition}
\label{prop:KO}
\cite[K\"uhn--Osthus]{KuOs06soda}
Let $r, n\ge 3$ be integers and $\varepsilon \in (0,1)$.
Let $G$ be an $n$-vertex graph with $\delta(G)\ge (1-\frac1r)n - n^{\varepsilon}$.
Then the decision problem for $K_r$-factor in $G$ is \texttt{NP-C}.
\end{proposition}

Therefore, 
Theorem~\ref{thm:main} bridges between this small gap and shows that the problem is \texttt{FPT} parameterized in $c$ (in particular, the problem is in \texttt{P} when $c=O(\log^{r^{-(r+4)}} n)$).

A natural further extension is the algorithmic question on $K_r$-tilings in graphs of arbitrary size, which has been studied by Fellows, Knauer, Nishimura, Ragde, Rosamond, Stege, Thilikos, and Whitesides~\cite{FellowsFPT} who showed that (among others) in general graphs the search/decision problem for a $K_r$-tiling of size $c$ is \texttt{FPT} parameterized in $c$.
\begin{theorem}
\cite[Theorem 2]{FellowsFPT}
\label{thm:fellows}
Given an $r$-vertex graph $F$ and $n$-vertex graph $G$ and integer $k > 0$, there is an algorithm that determines whether $G$ has an $F$-tiling of size $k$ in time $O(n^{r}+2^{O(kr)})$.
\end{theorem}

It will be interesting to establish a result similar to Theorem~\ref{thm:main} for $K_r$-tilings of arbitrary size.
Together with a simple reduction to the factor case, our result yields the following result on the mixture of tilings.

Given a graph $G$, an \emph{$\{mK_r, m'K_{r-1}\}$-factor} in $G$ consists of a $K_r$-tiling of size $m$ and a $K_{r-1}$-tiling of size $m'$ which are vertex-disjoint and cover all the vertices of $G$. 
Theorem \ref{thm:main} can be extended to decide whether a graph contains an $\{mK_r, m'K_{r-1}\}$-factor.

\begin{theorem}\label{thm:mixed tiling}
Given $r, c, n, m,m'\in \mathbb{N}$ such that $m\in [0, \frac{n}{r}]$, $m'\in [0, \frac{n}{r-1}]$ and $mr+m'(r-1)=n$.
Let $G$ be an $n$-vertex graph with $\d (G) \ge (1-\frac{1}{r-1}) n+\frac{m}{r-1} - c$.
Then there is an algorithm with running time $n^{O(1)}2^{O(c^{r^{r+4}})}$ that decides whether $G$ contains an $\{mK_r, m'K_{r-1}\}$-factor, that is, the problem is \texttt{FPT} parameterized in $c$.
Moreover, the algorithm either returns an $\{mK_r, m'K_{r-1}\}$-factor or a certificate showing that none exists.
\end{theorem}

Unfortunately, it does not solve the decision problem for $K_r$-tilings.
More specifically, $G$ may have a $K_r$-tiling of size $m$ even if it has no $\{mK_r, m'K_{r-1}\}$-factor.
Nevertheless, the complement version of Theorem~\ref{thm:mixed tiling} gives a nice consequence on equitable colorings (Corollary~\ref{cor:1.5}).

\begin{proof}[Proof of Theorem \ref{thm:mixed tiling} using Theorem~\ref{thm:main}]
Let $A$ be a set of $m'$ ``omni'' vertices disjoint with $V(G)$, and let $G'$ be obtained by connecting all vertices of $A$ to all vertices of $G$.
Note that the vertices of $A$ are independent.
Note that $G'$ contains a $K_r$-factor if and only if $G$ contains an $\{mK_r, m'K_{r-1}\}$-factor.
Indeed, if $G'$ contains a $K_r$-factor $K$, then we remove the $m'$ vertices of $A$ from $K$ and obtain an $\{mK_r, m'K_{r-1}\}$-factor of $G$ (using the fact that the vertices in $A$ are independent, so no two of them are in the same $K_r$ in $K$). If $G$ contains an $\{mK_r, m'K_{r-1}\}$-factor, then $G'$ contains a $K_r$-factor because every vertex of $A$ connects to all vertices of $G$.

For every vertex $v\in V(G)$ we have
\[d_{G'}(v)\ge \d (G)+m' \ge  (1-\frac{1}{r-1}) n+\frac{m}{r-1} - c +m' = (1-\frac 1r)(n+m') - c\]
and for every vertex of $A$ its degree is $n$ in $G'$. Since $n\ge(1-\frac 1r)(n+m') - c$, we have $\d (G') \ge(1-\frac 1r)(n+m') - c.$
By applying Theorem~\ref{thm:main} on $G'$, the algorithm with running time $(n+m')^{O(1)}2^{O(c^{r^{r+4}})}$ decides whether $G'$ contains a $K_r$-factor.
Hence the algorithm also decides whether $G$ contains an $\{mK_r, m'K_{r-1}\}$-factor.
\end{proof}

By considering the complement graph of $G$ and letting $\Delta:=m+m'$ in Theorem \ref{thm:mixed tiling}, we have the following result on equitable colorings. 

\begin{corollary}\label{cor:1.5}
Given $r, c, n, \Delta\in \mathbb{N}$ such that $\frac{n}{r}\le \Delta\le \frac{n}{r-1}$. Let $G$ be an $n$-vertex graph with $\Delta (G) \le \Delta+c$. Then there is an algorithm with running time $n^{O(1)}2^{O(c^{r^{r+4}})}$ that decides whether $G$ has an equitable $\Delta$-coloring, that is, the problem is \texttt{FPT} parameterized in $c$. Moreover, the algorithm either returns an equitable $\Delta$-coloring or a certificate showing that none exists.   
\end{corollary}


\subsection{Related work}
Fomin, Golovach, Sagunov, and Simonov~\cite{Fomin22} obtained a collection of results on the existence of long cycles in 2-connected graphs with minimum degree assumptions or its slight weakening, including the algorithmic extension of Dirac's Theorem mentioned earlier.
In a followup paper, 
they \cite{MR4481798} also showed such an algorithmic extension of Erd\H{o}s--Gallai Theorem, resolving an open problem.
That is, they showed that for $c\in \mathbb N$, the decision problem of whether the circumference of a graph $G$ is at least $d(G)+c$ is \texttt{FPT} parameterized in $c$, where $d(G)$ is the average degree of $G$.
See also~\cite{FGKSS23, fomin2023approximating, Fomin-2024soda} for further results.

Regarding Hamiltonicity, Dahlhaus, Hajnal, and Karpi\'nski~\cite{DHK} showed that it is \texttt{NP-C} to decide whether an $n$-vertex graph $G$ with $\delta(G) \ge (1/2-\varepsilon) n$ contains a Hamiltonian cycle.
For graphs $G$ with $\delta(G)\ge n/2-c$, Jansen, Kozma, and Nederlof~\cite{JKN} in 2019 gave an algorithm that runs in time $2^{O(c)} n^{O(1)}$ which decides if $G$ contains a Hamiltonian cycle.

Another contribution to Question~\ref{que} is made by Han and Keevash~\cite{han2022finding}, who studied Question~\ref{que} regarding the minimum $(k-1)$-degree assumption that forces matchings of certain size in $k$-uniform hypergraphs.
More precisely, they showed that for $c\in \mathbb N$, the decision problem of whether a $k$-uniform hypergraph $H$ containing a matching of size $\delta+c$ is \texttt{FPT} parameterized in $c$, where $\delta$ is the minimum $(k-1)$-degree of $H$.

On the other hand, there is an example showing that Question~\ref{que} does not always enjoy an \texttt{FPT} transition from \texttt{P} to \texttt{NP-C} (assuming \texttt{P$\neq$NP}).
Indeed, a $k$-uniform tight cycle is a $k$-uniform hypergraph whose vertices can be cyclically ordered such that the edges are precisely the sets of consecutive $k$ vertices.
For $k=3$, it is proved by R\"odl, Ruci\'nski, Szemer\'edi~\cite{RRS11} that for large enough $n$, every 3-uniform hypergraph with minimum 2-degree at least $\lfloor n/2\rfloor$ contains a tight Hamiltonian cycle; in contrast, Garbe and Mycroft~\cite{GaMy} showed that there exists a constant $C>0$, such that it is \texttt{NP-hard} to decide whether a 3-uniform hypergraph with minimum 2-degree at least $ n/2-C$ contains a tight Hamiltonian cycle.

\subsection{Proof ideas}

Our proof uses a number of tools.
For parameterized algorithms, we use the color-coding technique from Alon, Yuster and Zwick~\cite{AYZ95} for searching small $K_r$-tilings and the perfect hash functions from Fredman, Kolm\'os, and Szemer\'edi~\cite{FKS}.
To illustrate our proof strategy we first state the following result, which clearly implies Theorem~\ref{thm:main}.

\begin{lemma}
\label{lem:char}
Let $G$ be an $n$-vertex graph with $\delta(G)\ge (1-1/r)n-c$.
Then in time $n^{O(1)}2^{O(c^{{r}^{r+2}})}$ (resp. $n^{O(1)}2^{O(c)}$ for $r=3$) one can find a family $\mathcal F$ of $K_r$-tilings each of size $O(c^{r^{r+2}})$ (resp. $O(c)$ for $r=3$) defined on $V(G)$ such that $G$ has a $K_r$-factor if and only if there is $F\in \mathcal F$ such that $F\subseteq G$.
Moreover, the existence of such $F$ in $G$ can be checked in time $n^{O(1)}2^{O(c^{{r}^{r+2}})}$ (resp. $n^{O(1)}2^{O(c)}$ for $r=3$), and the algorithm outputs a $K_r$-factor if one exists.
\end{lemma}

Indeed, we shall show that one can take $\mathcal F$ as a certain family of $K_r$-tilings of size $O(c^{r^{r+2}})$ (resp. $O(c)$ for $r=3$), which can be searched by the color-coding method.
Then an immediate question will be how to make it a \textit{necessary} condition.
This is answered by our Lemma~\ref{lem:slack}.
Roughly speaking, it says that given such a rigid structure, one can conclude that if $G$ has a $K_r$-factor, then it has a small $K_r$-tiling of similar property.

\begin{lemma}
[Slack Lemma]
\label{lem:slack}
Let $r\geq r'\geq 2$, $I \subseteq [r']$.
Let $t_1,t_2,\ldots, t_{r'}$ be $r'$ integers such that $\sum_{i=1}^{r'}t_i=0$ and $\sum_{i\in I} |t_i|>4r-3$.
Let $G$ be an $n$-vertex graph and $(V_1, V_2, \dots, V_{r'})$ be a partition of $V(G)$ with $|V_i|=n/r+t_i$ for each $i\in [r'-1]$ and $|V_{r'}|=(r-r'+1)n/r+t_{r'}$.
Suppose that $G$ contains a $K_r$-factor $K$.
Denote $t_I = \sum_{i\in I} |t_i|$. 
Then we can find a subset $K' \subseteq K$ of size at most $g(t_I)$ such that there exists $x\in \mathbb N$ with $|V_i\cap V(K')|=x+t_i$, $i\in [r'-1]\cap I$ and $|V_{r'}\cap V(K')|=(r-r'+1)x+t_{r'}$ if $r'\in I$, where $g(t_I)=(3t_I)^r$ for $r\ge 4$ and $g(t_I)=8t_I$ for $r=3$.
In particular, $G-K'$ still contains a $K_r$-factor. 
\end{lemma}

Lemma~\ref{lem:slack} is proved by some careful analysis together with the pigeonhole principle.
In our proof, when we are given such a partition of $V(G)$, we always try to make the parts balanced, that is, the ideal partition should have part sizes with ratio $1:\cdots:1:(r-r'+1)$.
Then Lemma~\ref{lem:slack} says that $G$ has a $K_r$-factor only if it has a small $K_r$-tiling $K'$ that eliminates the discrepancy of coordinates in $I$ (if $I=[r']$, then it makes all parts balanced).

In our proof, we shall build such a vertex partition so that $t_I=t_I(c)=O(c^{r^{r+2}})$ (resp. $O(c)$ for $r=3$) and show that for almost all cases, this existence of $K'$ is also a sufficient condition for the existence of $K_r$-factor; in some rare cases, a slightly more complicated condition works.

Now let us elaborate more on our proof of sufficiency.
First, it is easy to see that an $n$-vertex graph with an independent set of size $n/r+1$ does not have a $K_r$-factor.
This simple construction serves as the extremal construction for Theorem~\ref{thm:HSze}.
A popular proof strategy in extremal graph theory aiming for exact result is the stability method, which exploits the structural information of the host graph.
Tailored to our problem, we indeed show that in Theorem~\ref{thm:main} if a graph $G$ is far from the extremal construction, namely, $G$ does not contain a sparse set of size $n/r$, then $G$ \emph{must} contain a $K_r$-factor.
The proof follows the stability method and the absorbing method.

For the remaining case ($G$ has at least one sparse $n/r$-set), it is crucial to know how many vertex-disjoint sparse $n/r$-sets $G$ contains.
For this we use Algorithm~\ref{alg:SparseSet} to find such a partition $(A_1,\dots,A_s,B)$ and edit it slightly.
Let us denote the resulting vertex partition of $G$ by $(A_1^1,\dots,A_s^1,B^1)$ (see its definition in Section 6), where $s$ is the maximum number of vertex-disjoint sparse $n/r$-sets $G$ has.
All the $A_i^1$ sets have size roughly $n/r$ (and thus $|B^1|\approx (r-s)n/r$).
So let us call a copy of $K_r$ balanced if it contains exactly one vertex from each $A_i^1$ and $r-s$ vertices from $B^1$.
Since each $A_i$ is sparse, we infer that each $G[A_i, A_j]$ and $G[A_i, B]$ is almost complete (and so do $G[A_i^1, A_j^1]$ and $G[A_i^1, B^1]$).
This motivates the ``ideal case'' of our proof (Lemma~\ref{lem: ideal}).
To reduce our $G$ to this ideal case, we need to (1) make the partition balanced, that is, its sizes have ratio $1:1:\cdots:1:(r-s)$ and (2) cover the vertices with wrong degree (in the sense of Lemma~\ref{lem: ideal}) by a small $K_r$-tiling with balanced copies.
However, (1) is not always possible and indeed it reduces to finding a small $K_r$-tiling of certain type and we shall solve it by the color-coding technique developed by Alon--Yuster--Zwick.
For a good summary on this strategy see our characterization results (Lemmas~\ref{lem:charact} and~\ref{lem:charact2}) in Section 6, which give the full version of Lemma~\ref{lem:char}.

\subsection{Sharpness of the result}

We build the first parameterized algorithm for $K_r$-\textsc{Factor In Graphs With High Minimum Degree} that bridges the gap from $c=0$ (\texttt{P}) to $c=n^{\Theta(1)}$ (\texttt{NP-C}).
However, the running time $n^{O(1)}2^{O(c^{{r}^{r+4}})}$ of our parameterized algorithm can probably be improved -- it will be interesting to design more efficient algorithms with reduced exponent of $c$, or even of type $n^{O(1)}2^{O(c)}$ (we have already achieved this type for $r=3$), which would yield a polynomial-time solution for $c=O(\log n)$.

\subsection*{Notation}
For positive integers $a,b$ with $a<b$, let $[a]=\{1,2,\ldots,a\}$ and $[a,b]=\{a,a+1,\ldots,b\}$.
For constants $x,y,z$, $x=y\pm z$ means that $y-z\le x\le y+z$. We write $x \ll y \ll z$ to mean that we choose constants from right to left, that is, for any $z>0$, there exist functions $f$ and $g$ such that, whenever $y \le f(z)$
and $x \le g(y)$, the subsequent statement holds. Hierarchies of other lengths are defined analogously.

For a graph $G$, we use $e(G)$ to denote the number of edges in $G$.
For a vertex $v\in V(G)$ and a vertex subset $A$, let $d_A(v)=|N_G(v)\cap A|$.
For two disjoint vertex subsets $A,B\subseteq V(G)$, let $e(A,B)=e(G[A,B])$.
Let $K$ be a $K_r$-tiling in $G$. Then we use $|K|$ to denote the size of $K$, that is, the number of vertex-disjoint copies of $K_r$ it contains.

\section{Non-extremal case}

In this section we study the case when the host graph $G$ has no large sparse set.
We use the absorbing method, pioneered by R\"odl, Ruci\'nski and Szemer\'edi~\cite{RRS06} in the early 2000s and used implicitly in the 1990s.
It usually reduces the task of building spanning subgraphs into building an absorbing set and an almost spanning subgraph.

\begin{definition}[Sparse set]
A vertex subset $U$ in an $n$-vertex graph $G$ is called a \emph{$\gamma$-sparse $k$-set} of $G$ if $|U|=k$ and $e(G[U])< \gamma n^2$.
\end{definition}

In this short section we prove the following result assuming two supporting lemmas.

\begin{theorem}[Non-extremal case]\label{thm:non-ext}
Suppose $r\ge 3$ and $1/n\ll \a\ll \g\ll 1/r$.
Let $G$ be an $n$-vertex graph with $\delta(G)\ge (1-1/r-\a)n$.
If $G$ has no $\gamma$-sparse $\frac{n}{r}$-set, then there is an algorithm running in time $O(n^{4r^2})$ that finds a $K_r$-factor in $G$.
\end{theorem}

To state our absorbing lemma, we first give the definitions of absorbers and absorbing sets.
Roughtly speaking, absorbers are constant-size sets that can contribute to the $K_r$-factor in multiple ways.

\begin{definition}[Absorber]
For any $r$-set $U\subseteq V(G)$, we say that a set $A_U\subseteq V(G)\setminus U$ is an \emph{absorber} for $U$ if $|A_U|= r^2$ and both $G[A_U]$ and $G[A_U\cup U]$ contain $K_r$-factors.
\end{definition}

An absorbing set is a vertex set satisfying that every $r$-set of vertices has many absorbers in it.

\begin{definition}[Absorbing set]
A vertex set $A\subseteq V(G)$ is called a $\xi$-\emph{absorbing set} for some $\xi>0$ if it consists of at most $\xi n$ vertex-disjoint $r^2$-sets and every $r$-set has at least $\xi^2 n$ absorbers in $A$.
\end{definition}


The following two lemmas, the absorbing lemma and the almost cover lemma, are by now regular routes for building spanning subgraphs.

\begin{lemma}[Absorbing lemma]\label{lem:abs}
Suppose $r\ge 3$ and $1/n \ll \xi \ll \a \ll \g \ll 1/r$.
If $G$ is an $n$-vertex graph with $\delta(G)\ge (1-1/r-\a)n$ and $G$ has no $\gamma$-sparse $\frac{n}{r}$-set,
then we can find a $\xi$-absorbing set $A$ in time $O(n^{4r^2})$.
\end{lemma}

\begin{lemma}[Almost cover]\label{lem:almost}
Suppose $r\ge 3$ and  $1/n \ll \tau\ll \a \ll \g \ll 1/r$.
If $G$ is an $n$-vertex graph with $\delta(G)\ge (1-1/r-\a)n$ and $G$ has no $\gamma$-sparse $\frac{n}{r}$-set, then we can find in time $O(n^{r+1})$ a $K_r$-tiling covering all but at most $\tau n$ vertices of $G$.
\end{lemma}

We first use the two lemmas to establish the non-extremal case.

\begin{proof}[Proof of Theorem \ref{thm:non-ext}]
Choose constants \[1/n\ll \tau\ll \xi\ll \a \ll \g\ll 1/r.\]

We first apply Lemma \ref{lem:abs} on $G$ and obtain in time $O(n^{4r^2})$ a $\xi$-absorbing set $A$ consisting of at most $\xi n$ vertex-disjoint  $r^2$-sets such that every $r$-set has at least $\xi^2 n$ absorbers in $A$.
Let $n'=|V(G)\setminus A|$. Since $|A|\le \xi r^2 n$ and $\xi\ll \a$, we have  $n'\ge (1-\xi r^2)n$ and $\delta(G- A)\ge (1-1/r-2\a)n\ge (1-1/r-2\a)n'$.
For an $\frac{n'}{r}$-set in $V(G)\setminus A$, it induces at least $\g n^2-\xi rn^2\ge \g n^2/2$ edges and so $G- A$ has no $\frac{\gamma}{2}$-sparse $\frac{n'}{r}$-set.
Then we apply Lemma \ref{lem:almost} to $G- A$ with $2\alpha$ in place of $\alpha$ and $\g/2$ in place of $\g$, and obtain in time $O(n^{r+1})$ a $K_r$-tiling $\mathcal{K}$ covering all but at most $\tau n'$ vertices of $G- A$. Denote by $U$ the set of uncovered vertices in $G- A$.
We arbitrarily divide $U$ into at most $\tau n'/r$ vertex-disjoint $r$-sets.
Since $\tau\ll \xi^2$, we can greedily assign each $r$-set a unique absorber in $A$. Hence $G[A\cup U]$ contains a $K_r$-factor. Together with $\mathcal{K}$, we obtain a $K_r$-factor of $G$.
\end{proof}

We prove the two lemmas in the next two sections.

\section{Absorbing lemma}

In this section we prove Lemma \ref{lem:abs}.
The original proofs of the absorbing method build the desired absorbing set by probabilistic arguments, which yields randomized algorithms that output the absorbing set with high probability.
Following previous works, we use the following result to ``find'' the absorbing set, which is a typical application of the \emph{conditional expectation method} for derandomization.

\begin{lemma}\cite[Proposition 4.7]{GaMy} \label{lem:derandomization}
Fix constants $\beta_1>\beta_2>0$ and integers $m, M, N$ and $\ell\le N$ such that $N$ and $\ell$ are sufficiently large, and that $M\le (1/8)\exp(\beta_2^2\ell/(3\beta_1))$. Let $U$ and $W$ be disjoint sets of sizes $|U| = M$ and $|W| = N$. Let $H$ be a graph with vertex set $U\cup W$ such that $H[U]$ is empty, $H[W]$ has exactly $m$ edges, and $d_H(u)\ge \beta_1 N$ for every $u\in U$. Then we can find an independent set $L\subseteq W$ in $H$ in time $O(N^4+MN^3)$ such that $(1-\nu)\ell\le |L|\le \ell$ and $|N_H(u)\cap L|\ge (\beta_1 - \beta_2 -\nu)\ell$ for every $u\in U$, where $\nu = 2m\ell/N^2$.
\end{lemma}

In applications, the edges between $U$ and $W$ are the useful edges, which are usually used to model the desired property, e.g., every $u\in U$ ``has many absorbers'' in $W$; and the edges in $W$ model the relation that we would like to avoid.

\begin{proof}
[Proof of Lemma \ref{lem:abs}]
Fix an integer $r\ge 3$ and choose constants $1/n\ll \xi\ll \a\ll \g\ll 1/r$.

The proof follows the standard proof of the absorbing method, starting with the enumeration of the absorbers.
We claim that any $r$-set $\{v_1,\ldots,v_r\}\subset V(G)$ has $\Omega(n^{r^2})$ absorbers.
Indeed, we first choose $r$ distinct vertices $u_1,\ldots,u_r\in V(G)\setminus \{v_1,\ldots,v_r\}$ one by one such that $\{u_1,\ldots,u_r\}$ induces a copy of $K_r$.
Since any $r-1$ vertices have at least $(r-1)(1-1/r-\a)n-(r-2)n\ge(1/r-\a r)n$ common neighbors in $G$, there are at least $(1/r-\a r)^r n^r/r!$ choices for $u_1,\ldots,u_r$.
Then for each $i\in [r]$ we choose $r-3$ distinct vertices $w_{i1},w_{i2},\ldots,w_{i,{r-3}}$ from the common neighborhood of $v_i$ and $u_i$ such that both $v_i,w_{i1},w_{i2},\ldots,w_{i,{r-3}}$ and $u_i,w_{i1},w_{i2},\ldots,w_{i,{r-3}}$ induce copies of $K_{r-2}$.
Finally, we choose an edge $x_iy_i$ from the common neighborhood $U_i$ of $v_i,u_i,w_{i1},w_{i2},\ldots,w_{i,{r-3}}$ for each $i\in [r]$.
Since $|U_i|\ge (1/r-\a r)n$ and $G$ has no $\gamma$-sparse $\frac{n}{r}$-set, $G[U_i]$ induces at least $(\g - \a r)n^2$ edges for any $i\in [r]$.
Hence $\bigcup_{i=1}^r\{u_i,w_{i1},w_{i2},\ldots,w_{i,{r-3}},x_i,y_i\}$ is an absorber for the $r$-set $\{v_1,\ldots,v_r\}$ and it has at least
\[(1/r-\a r)^r n^r \cdot \prod_{i=1}^r(1/r-\a r)^{r-3}n^{r-3}\cdot \prod_{i=1}^r(\g - \a r)n^2/r^2!\ge \gamma'n^{r^2}\]
choices, where $\gamma'= (1/r-\a r)^{r^2-2r}(\g - \a r)^r/r^2!$. Thus, every $r$-set in $V(G)$ has at least  $\gamma'n^{r^2}$ absorbers, which can be found (by brute force) in time $O(n^{r^2})$.


To apply Lemma \ref{lem:derandomization}, we define an auxiliary graph $H$ with $U={V(G)\choose r}$ 
and $W={V(G)\choose {r^2}}$ such that $u\in U$ and $w\in W$ are adjacent if and only if $w$ is an absorber for $u$, and $w_1,w_2\in W$ are adjacent if and only if $w_1\cap w_2\neq \emptyset$.
In the notation of Lemma \ref{lem:derandomization} we have $M={n\choose r}$, $N={n \choose {r^2}}$ and
\[m=e(H[W])\le {n \choose {r^2}}\cdot r^2 \cdot {n \choose {r^2-1}}=\frac{r^4 N^2}{n-r^2+1}.\]
Let $\beta_1=\g'$, $\beta_2=\g'/3$ and $\ell=\xi n$.
Since every $r$-set in $V(G)$ has at least  $\gamma'n^{r^2}$ absorbers, we get $d_H(u)\ge \beta_1 N$ for every $u\in U$ and
\[\exp\left(\frac{\beta_2^2\ell}{3\beta_1}\right)=\exp\left(\frac{\g'\xi n}{27}\right)\ge 8{n\choose r}=8M,\]
as $n$ is large enough.
By Lemma \ref{lem:derandomization}, in time $O(N^4+MN^3)=O(n^{4r^2}+n^rn^{3r^2})=O(n^{4r^2})$, we can find an independent set $L\subseteq W$ in $H$ such that $(1-\nu)\ell\le |L| \le \ell$ and $|N_H(u)\cap L|\ge (\beta_1 - \beta_2 -\nu)\ell$ for all $u\in U$, where
\[\nu = \frac{2m\ell}{N^2}\le \frac{2r^4 N^2\xi n}{(n-r^2+1)N^2}<\frac{\g'}{3},\]
as $\xi \ll \a \ll \g\ll 1/r$.
Note that $L$ is a family of vertex-disjoint $r^2$-sets of $V(G)$ by the definition of $H[W]$ and $|L|\le \xi n$.
Moreover, every $r$-set of $V(G)$ has at least
\[(\beta_1 - \beta_2 -\nu)\ell>\frac{\g'\xi n}{3}\ge \xi^2 n\] absorbers in $L$, as $\xi \ll \a \ll \g$.
Let $A\subseteq V(G)$ be the vertex set consisting of vertex-disjoint $r^2$-sets in $L$. Then $A$ is the desired $\xi$-absorbing set.
\end{proof}

\section{Almost cover}

The main tools for the proof of Lemma~\ref{lem:almost} are the regularity method and a result of Keevash and Mycroft on almost perfect matchings in clique-complexes.

\subsection{Regularity}
We use the standard graph regularity method. We first give related definitions and lemmas.
\begin{definition}[Regular pair]
Given a graph $G$ and disjoint vertex subsets $X,Y\subseteq V(G)$, the \emph{density} of the pair $(X,Y)$ is defined as \[d(X,Y):=\frac{e(X,Y)}{|X||Y|}.\] For $\e>0$, the pair $(X,Y)$ is $\e$-\emph{regular} if for any $A\subseteq X, B\subseteq Y$ with $|A|\ge \e |X|, |B|\ge \e |Y|$, we have \[|d(A,B)-d(X,Y)|<\e.\] Additionally, if $d(X,Y)\ge d$ for some $d\ge0$, then we say that $(X,Y)$ is $(\e,d)$-\emph{regular}.
\end{definition}
\begin{lemma}\cite[Slicing lemma]{kom-sim}\label{lem:slicing}
Let $(X,Y)$ be an $\e$-regular pair with density $d$, and for some $\eta>\e$, let $X'\subseteq X, Y'\subseteq Y$ with $|X'|\ge \eta |X|, |Y'|\ge \eta |Y|$. Then $(X',Y')$ is an $\e'$-regular pair with $\e'=\max\{\e/\eta,2\e\}$, and for its density $d'$ we have $|d'-d|<\e$.
\end{lemma}

\begin{lemma}\cite[Counting lemma]{kom-sim} \label{lem:counting}
Given $d>\e>0$, $m\in \mathbb{N}$ and $H$ some fixed graph on $r$ vertices, let $G$ be a graph obtained by replacing every vertex $x_i$ of $H$ with an independent set $V_i$ of size $m$ and every edge of $H$ with an $(\e,d)$-regular pair on the corresponding sets. If $\e\le \frac{d^r}{(2+r)2^r}=:d_0$, then there are at least $(d_0m)^r$ embeddings of $H$ in $G$ so that each $x_i$ is embedded into the set $V_i$.
\end{lemma}

\begin{lemma}\cite[Degree form of the Regularity Lemma]{algorithmicRegLem,kom-sim}\label{lem:reg}
For every $\e>0$, there is an $M=M(\e)$ such that if $G=(V,E)$ is an $n$-vertex graph and $d\in[0,1]$ is any real number, then we can find in time $O(n^{2.376})$ a partition $V=V_0\cup V_1\cup \cdots \cup V_k$ and a spanning subgraph $G'\subseteq G$ with the following properties:
\begin{itemize}
  \item $1/\e \le k\le M$,
  \item $|V_0|\le \e |V|$ and $|V_i|=m$ for all $1\le i\le k$ with $m\le \e |V|$,
  \item $d_{G'}(v)>d_{G}(v)-(d+\e)|V|$ for all $v\in V$,
  \item $e(G'[V_i])=0$ for all $i\ge 1$,
  \item all pairs $(V_i,V_j)$ $(1\le i<j\le k)$ are $\e$-regular in $G'$ with density $0$ or at least $d$, where in the latter case, we indeed have $E_{G'}(V_i, V_j)=E_{G}(V_i, V_j)$.
\end{itemize}
\end{lemma}
\noindent In applications, we usually call $V_0,V_1,\ldots,V_k$ \emph{clusters} and call the cluster $V_0$ \emph{exceptional set}.

\begin{definition}[Reduced graph]
Given an arbitrary graph $G=(V,E)$, a partition $V=V_1\cup \cdots \cup V_k$, and reals $\e,d>0$, the \emph{reduced graph} $R=R(\e,d)$ of $G$ is defined as follows:
\begin{itemize}
  \item $V(R)=\{V_1,\ldots, V_k\}$,
  \item $V_iV_j\in E(R)$ if and only if $(V_i,V_j)$ is $(\e,d)$-regular.
\end{itemize}
\end{definition}

When applying the Regularity Lemma, we start with an $n$-vertex graph $G=(V,E)$ and parameters $\e,d>0$, and then obtains a partition $V=V_0\cup V_1\cup \cdots \cup V_k$ and a subgraph $G'$ as in Lemma \ref{lem:reg}. Then we usually study the properties of reduced graph $R=R(\e,d)$ of $G'\setminus V_0$. By Lemma \ref{lem:reg},
\[
\delta(R)\ge \frac{\left(\delta(G)-(d+\e)|V|-|V_0|\right)m}{m^2} \ge \frac{\delta(G)-(d+2\e)n}{m}.
\]
In particular, if $\delta(G)\ge cn$, then $\delta(R)\ge (c-d-2\e)k$.
Moreover, we show that $R$ contains a sparse set if and only if $G$ contains a sparse set.
Recall that a vertex subset $U$ in an $n$-vertex graph $G$ is called a $\gamma$-sparse $k$-set of $G$ if $|U|=k$ and $e(G[U])< \gamma n^2$.

\begin{fact}\label{lem:sparse set1}
Let $0<\beta<1$ be a real number. If $R$ has a $\beta$-sparse $\frac{k}{r}$-set, then $G$ has a $\left(\beta+2\e+d\right)$-sparse $\frac{n}{r}$-set.
\end{fact}

\begin{proof}
Let $U$ be a $\beta$-sparse $\frac{k}{r}$-set of $R$ and we restore $U$ to a sparse set of $G$.
Let $W$ be the union of the clusters in $U$.
Since $|V_i|=m\ge \frac{(1-\e)n}{k}$ for every $i\in[k]$, we have $|W|=m|U|=\frac{mk}{r}\ge \frac{(1-\e)n}{r}$.
Since $e(R[U])\le \beta k^2$ and $d_{G'}(v)>d_{G}(v)-(d+\e)n$ for all $v\in V$, we have
\begin{align}
e(G[W])&\le \beta k^2m^2+ (e(G) - e(G'))\nonumber\\
       &\le \beta n^2+(d+\e)n^2. \nonumber
\end{align}
We add up to $\frac{\e n}{r}$ vertices of $V(G)\setminus W$ to $W$ and obtain an $\frac{n}{r}$-set $W'$. Then $|W'|=\frac{n}{r}$ and $e(G[W'])\le \beta n^2+(d+\e)n^2+\frac{\e n^2}{r^2}\le (\beta+2\e+d) n^2$. That is, $W'$ is a $(\beta+2\e+d)$-sparse $\frac{n}{r}$-set of $G$.
\end{proof}

\begin{fact}\label{lem:sparse set2}
Let $0<\eta<1$. If $R$ has no $\left(\frac{8\eta}{\e^2d}+\frac{3\e}{r}\right)$-sparse $\frac{k}{r}$-set, then $G$ has no $\eta$-sparse $\frac{n}{r}$-set.
\end{fact}

\begin{proof}
Suppose not. Let $W$ be an $\eta$-sparse $\frac{n}{r}$-set of $G$.
Let $W_i=W\cap V_i$ for each $i\in [k]$ and $U=\{V_i\in V(R): |W_i|>\e m\}$. Then
\[
|U|\ge \frac{|W|-\e mk-|V_0|}{m}\ge \frac{k}{r}-\e k-\frac{\e n}{m}\ge \frac{k}{r}-\e k-\frac{\e k}{1-\e}\ge \frac{k}{r}-3\e k.
\]
Now we calculate $e(R[U])$. Note that if $V_iV_j$ is an edge of $R[U]$, then $(V_i,V_j)$ is $\e$-regular with density $d_{ij}\ge d$. Hence we have $\sum_{V_iV_j\in E(R[U])}d_{ij}|V_i||V_j|\ge e(R[U])\cdot dm^2$. Moreover, since $|W_i|,|W_j|>\e m$, by Lemma \ref{lem:slicing}, $(W_i,W_j)$ is regular with density $d_{ij}'$, where $d_{ij}'=(1\pm \e)d_{ij}$. Then
\[
\sum_{V_iV_j\in E(R[U])}d_{ij}|V_i||V_j|
   < \sum_{V_iV_j\in E(R[U])}\frac{d_{ij}|W_i||W_j|}{\e^2}
   \le \sum_{V_iV_j\in E(R[U])}\frac{e_G(W_i,W_j)}{\e^2(1- \e)}\le \frac{\eta n^2}{\e^2(1-\e)}.
\]
Hence
\[e(R[U])\le \frac{\eta n^2}{\e^2(1-\e)dm^2}\le \frac{\eta}{\e^2(1-\e)d}\cdot\frac{k^2}{(1-\e)^2}<\frac{8\eta k^2}{\e^2d}.\]
We add up to $3\e k$ vertices of $V(R)\setminus U$ to $U$ and obtain a $\frac{k}{r}$-set $U'$. Then $|U'|=\frac{k}{r}$ and $e(R[U'])\le \frac{8\eta k^2}{\e^2d}+3\e k\cdot\frac{k}{r}=\left(\frac{8\eta}{\e^2d}+\frac{3\e}{r}\right)k^2$. That is, $U'$ is a $\left(\frac{8\eta}{\e^2d}+\frac{3\e}{r}\right)$-sparse $\frac{k}{r}$-set of $R$, a contradiction.
\end{proof}

\subsection{Almost perfect matching in $r$-complexes}
In the proof of Lemma \ref{lem:almost}, we will use the result in \cite{KM1} about the existence of an almost perfect matching in an $r$-complex. To state the results, we first introduce some definitions. An \emph{$r$-system} is a hypergraph $J$ in which every edge of $J$ has size at most $r$ and $\emptyset \in J$.
An edge of size $k$ ($k \le r$) in $J$ is called a \emph{$k$-edge} of $J$. We use $J_k$ to denote the $k$-uniform hypergraph on $V(J)$ with all $k$-edges of $J$.

\begin{definition}[$r$-complex]
An \emph{$r$-complex} $J$ is an $r$-system whose edge set is closed under inclusion, that is, if $e\in J$ and $e'\subseteq e$ then $e'\in J$.
\end{definition}

In an $r$-complex $J$, the \emph{degree} $d(e)$ of an edge $e$ is the number of $(|e|+1)$-edges $e'$ of $J$ which contain $e$ as a subset. The \emph{minimum $k$-degree} of $J$, denoted by $\delta_k(J)$, is the minimum of $d(e)$ taken over all $k$-edges $e\in J$. Since the edge set of $J$ is closed under inclusion, we have $\delta_k(J)\le \delta_{k-1}(J)$ for each $k\in [r-1]$. The \emph{degree sequence} of $J$ is
\[\delta(J)=(\delta_0(J),\delta_1(J),\ldots,\delta_{r-1}(J)).\]

\begin{definition}[Clique $r$-complex]
The \emph{clique $r$-complex} $J$ of a graph $G$ is defined by taking $J_k$ to consist of the $k$-cliques of $G$ for $0\le k\le r$.
\end{definition}

\begin{construction}[Space barrier]
Suppose $V$ is a set of $n$ vertices, $j\in [r-1]$ and $S\subseteq V$. Let $J(S,j)$ be the $r$-complex in which $J(S,j)_i$ (for $0\le i\le r$) consists of all $i$-sets in $V$ that contain at most $j$ vertices of $S$.
\end{construction}

\begin{definition}
Suppose $G$ and $H$ are $r$-uniform graphs on the same set $V$ of $n$ vertices and $0<\beta<1$. We say that $G$ is \emph{$\beta$-contained} in $H$ if all but at most $\beta n^r$ edges of $G$ are edges of $H$.
\end{definition}

The following theorem is proved by Keevash and Mycroft~\cite{KM1}.

\begin{theorem}\cite{KM1}\label{thm:complex}
Suppose that $1/n\ll 1/\ell\ll \a\ll \beta\ll 1/r$. Let $J$ be an $r$-complex on $n$ vertices such that
\[\delta(J)\ge \left(n,\left(1-\frac{1}{r}-\a\right)n,\left(1-\frac{2}{r}-\a\right)n,\ldots,\left(\frac{1}{r}-\a\right)n\right).\]
Then $J$ has at least one of the following properties:
\begin{enumerate}
  \item $J_r$ contains a matching that covers all but at most $\ell$ vertices.
  \item $J_r$ is $\beta$-contained in $J(S,j)_r$ for some $j\in [r-1]$ and $S\subseteq V(J)$ with $|S|=\lfloor jn/r\rfloor$.
\end{enumerate}
\end{theorem}

We will give a corollary of Theorem~\ref{thm:complex}, i.e. Corollary~\ref{cor:complex}, which is easier to apply than Theorem~\ref{thm:complex}. The proof of Corollary~\ref{cor:complex} also use the following fact.

\begin{fact}\label{fact: counting Kr in graph with no sparse set}
For any $\gamma_H\in (0,1)$, there exist $\varepsilon,\delta>0$ such that the following holds. Let $H$ be an $h$-vertex graph. If $e(H)\ge (1-\frac{1}{k-1}-\varepsilon)\frac{h^2}{2}$ and $H$ has no $\gamma_H$-sparse $\frac{h}{k-1}$-set, that is, every $\frac{h}{k-1}$-set in $V(H)$ induces at least $\gamma_H h^2$ edges, then $H$ has at least $\delta h^k$ copies of $K_k$.
\end{fact}

To prove Fact~\ref{fact: counting Kr in graph with no sparse set}, we shall use the graph removal lemma by Erd\H{o}s, Frankl and R\"odl~\cite{EFR86} and the Erd\H{o}s–Simonovits stability theorem~\cite{ES66}.
   
\begin{lemma}\cite{EFR86}\label{lem:removal}
    For any graph $ F $ and any $ \epsilon > 0 $, there exists $ \alpha > 0 $ such that any graph on $ n $ vertices which contains at most $ \alpha n^{v(F)} $ copies of $ F $ may be made $ F $-free by removing at most $ \epsilon n^2 $ edges.
\end{lemma} 
\medskip

 A classical result in extremal graph theory is the Erd\H{o}s–Simonovits stability theorem~\cite{ES66}, which says that for every $\epsilon>0$ there exists $\alpha>0$ such that if $G$ is a $K_{r+1}$-free graph on $n$ vertices with $e(G)\ge \mathrm{ex}(n,K_{r+1})-\alpha n^2$, then $G$ can be obtained from the Tur\'{a}n graph $T_k(n)$ by adding and deleting at most $\epsilon n^2$ edges. Now we are ready to prove Fact~\ref{fact: counting Kr in graph with no sparse set}.

\begin{proof}
[Proof of Fact~\ref{fact: counting Kr in graph with no sparse set}]
Given $\gamma_H\in (0,1)$, we choose $\delta,\varepsilon\ll\varepsilon_1\ll\varepsilon_2\ll \gamma_H$.
Suppose that $H$ has less than $\delta h^k$ copies of $K_k$, by Lemma~\ref{lem:removal} with $\epsilon=\varepsilon_1$ and $\alpha=\delta$, $H$ can be made $K_{k}$-free by removing at most $\varepsilon_1 h^2$ edges. 
By Erd\H{o}s--Simonovits stability theorem~\cite{ES66} with $\epsilon=\varepsilon_2$ and $\alpha=\varepsilon/2$, the remaining graph can be changed to the Tur\'{a}n graph $T_{k-1}(h)$ by altering at most $\varepsilon_2h^2$ adjacencies. Hence $H$ has a $\gamma_H$-sparse $\frac{h}{k-1}$-set, a contradiction.
\end{proof}

Now we state the corollary of Theorem~\ref{thm:complex} and give its proof.

\begin{corollary}\label{cor:complex}
Suppose that $1/n\ll 1/\ell\ll \a\ll \mu\ll 1/r$.
Let $G$ be an $n$-vertex graph with $\delta(G)\ge (1-1/r-\a)n$. Then either $G$ contains a $K_r$-tiling that covers all but at most $\ell$ vertices, or it contains a $\mu$-sparse $\lfloor\frac{n}{r}\rfloor$-set.
\end{corollary}
\begin{proof}
Choose constants $1/n\ll {1/\ell}\ll \a\ll \beta\ll \mu\ll 1/r$.
We consider the clique $r$-complex $J$ of $G$, that is, $J_i$ consists of the $i$-cliques of $G$ for $0\le i\le r$. Note that $J_2=G$. Since $\delta(G)\ge (1-1/r-\a)n$, we have
\[\delta(J)\ge \left(n,\left(1-\frac{1}{r}-\a\right)n,\left(1-\frac{2}{r}-2\a\right)n,\ldots,\left(\frac{1}{r}-(r-1)\a\right)n\right).\]
By Theorem \ref{thm:complex}, $J$ satisfies the property (1) or (2).
If $J$ satisfies (1), then $G$ contains a $K_r$-tiling covering all but at most $\ell$ vertices and we are done.

Now suppose $J$ satisfies (2), that is, $J_r$ is $\beta$-contained in $J(S,j)_r$ for some $j\in [r-1]$ and $S\subseteq V(J)$ with $|S|=\lfloor jn/r\rfloor$, then $J_r$ has at most $\beta n^r$ edges $e$ with $|e\cap S|\ge j+1$.
So $G$ has at most $\beta n^r$ copies of $K_r$ each of which contains at least $j+1$ vertices of $S$.
Note that a copy of $K_{j+1}$ in $G[S]$ is contained in at least
\[\frac{1}{(r-j-1)!}\prod_{i=j+1}^{r-1}\left(1-\frac{i}{r}-2\a i\right)n\ge \left(\frac{n}{2r^2}\right)^{r-j-1}\] copies of $K_r$ in $G$.
By double counting, $G[S]$ has at most
\[{r\choose {j+1}}\beta n^r\bigg/\left(\frac{n}{2r^2}\right)^{r-j-1}\le (2r^2)^r\beta n^{j+1}\]
copies of $K_{j+1}$.
Since $|S|=\lfloor jn/r\rfloor$, we have $\delta(G[S])\ge (1-1/r-\a)n-(n-|S|)=(1-1/j)|S|-\a n$ and so\[e(G[S])\ge (1-1/j)|S|^2/2-\a n|S|/2.\]
By Fact \ref{fact: counting Kr in graph with no sparse set}, $G[S]$ has a $\mu$-sparse $\frac{|S|}{j}$-set, i.e. a $\mu$-sparse $\lfloor\frac{n}{r}\rfloor$-set of $G$.
\end{proof}

\subsection{Proof of Lemma \ref{lem:almost}}
Choose constants \[1/n\ll\e\ll d\ll \tau , {1/\ell}\ll \a\ll\mu\ll \g\ll 1/r.\]

We first apply Lemma \ref{lem:reg} on $G$ with parameters $\e$ and $d$ to obtain a partition $V(G)=V_0\cup V_1\cup\cdots \cup V_k$ in time $O(n^{2.376})$, for some $1/\e \le k\le M$ with $|V_i|=m\ge \frac{(1-\e)n}{k}$ for every $i\in[k]$. Let $R=R(\e,d)$ be the reduced graph for this partition. Since $\e\ll d \ll\a$, we have $\delta(R)\ge (1-1/r-\a-d-2\e)k\ge (1-1/r-2\a)k$. Moreover, $R$ has no $\mu$-sparse $\lfloor\frac{k}{r}\rfloor$-set. Indeed, otherwise by Lemma \ref{lem:sparse set1}, $G$ has a $(\mu+2\e+d)$-sparse $\frac{n}{r}$-set. Since $\e\ll d \ll\mu\ll \g$, $G$ has a $\g$-sparse $\frac{n}{r}$-set, a contradiction.

Then we apply Corollary \ref{cor:complex} on $R$ and obtain that $R$ contains a $K_r$-tiling $\mathcal{K}$ that covers all but at most $\ell$ vertices, which can be built
in constant time (as $R$ has a finite order).
Let $K$ be a copy of $K_r$ in $\mathcal{K}$.
Without loss of generality, let $V(K)=\{V_1,V_2,\ldots, V_r\}$.
By Lemmas \ref{lem:slicing} and \ref{lem:counting}, for each $i\in [r]$ we can greedily embed in $G[V_1,V_2,\ldots, V_r]$ vertex-disjoint copies of $K_r$ that cover all but at most $\e m$ vertices of $V_i$ in time $O(m^{r+1})$.
For each copy of $K_r$ in $\mathcal{K}$, we do the same operation, which gives us in time $O(n^{r+1})$ a $K_r$-tiling covering all but at most
\[\e m (|V(R)|-\ell)+\ell m+|V_0|\le \e m(k-\ell)+\ell m+\e n \le 2\e n+\e\ell n\le \tau n\]
vertices of $V(G)$.

\section{Preparation for the extremal case}
We give and show several lemmas and algorithms in this section that will be used in the proof of the extremal case, that is, when $G$ contains large sparse sets.

The following definitions are key to our proofs, and motivated by earlier works on problems of similar flavor.
The idea behind it is that when $G$ contains large sparse sets, one can find and analyze its (rough) structure.
All vertex partitions considered in this paper have an implicit order on their parts.

Let $r\geq r'\geq 2$ and $\varepsilon <  1/r$.
Let $t_1,t_2,\ldots, t_{r'}$ be $r'$ integers such that $|t_i|\le \varepsilon n$ for each $i\in [r']$ and $\sum_{i=1}^{r'}t_i=0$.
Let $G$ be an $n$-vertex graph and $(V_1, V_2, \dots, V_{r'})$ be a partition of $V(G)$ with $|V_i|=n/r+t_i$ for each $i\in [r'-1]$ and $|V_{r'}|=(r-r'+1)n/r+t_{r'}$.
Let $\boldsymbol{b}=(1,\ldots,1,r-r'+1)$ be an $r'$-dimensional vector.
For any $i\in [r']$, let $\boldsymbol{u}_i$ be the $i$-th unit vector, i.e. $\boldsymbol{u}_i$ has 1 on the $i$-th coordinate and 0 on the other coordinates.
Let $F$ be a $K_r$-tiling of size $f$ in $G$.

\begin{definition}[Index vector]\label{def:index_vector}
The \emph{index vector} of $F$ is $\textbf{i}(F)=(x_1,x_2,\ldots,x_{r'})$ by choosing $x_i=|V(F)\cap V_i|$ for every $i\in [r']$.
A copy of $K_r$ in $G$ is \emph{balanced} if its index vector is $\textbf{b}$. 
\end{definition}

Compared with a balanced copy of $K_r$, an \emph{$(i,j)$-copy of $K_r$} has index vector $ \boldsymbol{b}+\boldsymbol{u}_i-\boldsymbol{u}_j$, where $i, j\in [r']$ are distinct. 
For a copy of $K_r$, we refer to its index vector as its \emph{type}. Note that there are $y:={{r+r'-1}\choose {r'-1}}$ different types of $K_r$.

\begin{definition}[Slack]\label{def:slack}
The \emph{slack} of $F$ is $\textbf{s}(F) = (s_1, s_2, \dots, s_{r'}):=\textbf{i}(F)-f\textbf{b}$. Note that $\sum_{i=1}^{r'}s_i=0$ and $\sum_{s_i>0}s_i\le (r-1)f$.
The \emph{slack} of the partition $(V_1, V_2, \dots, V_{r'})$ is $(t_1, t_2, \dots, t_{r'})=(|V_1|, |V_2|, \dots, |V_{r'}|)-\frac{n}{r}\textbf{b}$. 
      The partition $(V_1, V_2, \dots, V_{r'})$ is \emph{balanced} if $|V_1|=\cdots= |V_{r'-1}|= |V_{r'}|/{(r-r'+1)}$.
\end{definition}

\medskip
The definitions are natural descriptions of the vertex distributions of a $K_r$-tiling.
The key idea here is that the slack of a $K_r$-factor (which covers $V(G)=V_1\cup V_2\cup \cdots \cup V_{r'}$) is equal to the slack of the partition, each coordinate of which is small.
Then, a necessary condition for having a $K_r$-factor is to have a $K_r$-tiling with the same slack.
The question is how small we can take such a $K_r$-tiling, which is answered by  Lemma~\ref{lem:slack}. Below we restate it here (using slack) for convenience and give its proof.
Before proceeding, we introduce a piece of notation.
Let $I$ be a subset of $[r']$ and $\boldsymbol{v}$ be an $r'$-dimensional vector. We use $\boldsymbol{v}|_I$ to denote the $|I|$-dimensional vector obtained by restricting $\boldsymbol{v}$ to the coordinates in $I$.



\medskip

\noindent\textbf{Lemma 1.8.}
\emph{
Let $r\geq r'\geq 2$, $I \subseteq [r']$.
Let $t_1,t_2,\ldots, t_{r'}$ be $r'$ integers such that $\sum_{i=1}^{r'}t_i=0$ and $\sum_{i\in I} |t_i|>4r-3$.
Let $G$ be an $n$-vertex graph and $(V_1, V_2, \dots, V_{r'})$ be a partition of $V(G)$ with $|V_i|=n/r+t_i$ for each $i\in [r'-1]$ and $|V_{r'}|=(r-r'+1)n/r+t_{r'}$.
Suppose that $G$ contains a $K_r$-factor $K$
whose slack is $\textbf{s}(K) = (t_1, t_2, \dots, t_{r'})$.
Denote $t_I = \sum_{i\in I} |t_i|$. 
Then we can find a subset $K' \subseteq K$ of size at most $g(t_I)$ such that $\textbf{s}(K')|_I  = \textbf{s}(K)|_I $, where $g(t_I)= (3t_I)^r$ for $r\ge 4$ and $g(t_I)= 8t_I$ for $r=3$. In particular, $G-K'$ still contains a $K_r$-factor.}

\medskip

Note that in this lemma we require the conclusion to hold for any given $I$, and that the $K'$ only adjusts the slack of the parts corresponding to $I$. Indeed, $K'$ will be searched by the color-coding technique (Lemma~\ref{lem:color coding}).
If this lemma only holds for $I=[r']$, then the slack of $K'$ equals to the slack of the whole partition. However, the running time required to search for $K'$ might become larger than expected, that is, no longer being polynomial. 
In the proof of the extremal case, we first classify the parts by the coordinates of the slack, and let $I\subseteq [r']$ be the set of indices with \textit{small} (to be quantified) slack coordinates.
In order to find a $K_r$-tiling $H$ whose slack is equal to the slack of the partition, we first use color-coding to find a $K_r$-tiling $K'$ that adjusts the slack of the parts corresponding to $I$. If no such structure can be found, then $G$ contains no $K_r$-factor by the contrapositive of Lemma~\ref{lem:slack}. Otherwise, we will proceed to construct the remaining part of $H$. This remaining part is guaranteed to exist and can be found in polynomial time.

\begin{algorithm}
        \caption{$K_r$-tiling of small size}
        \label{alg:slack}
        \KwData{a $K_r$-factor $K$}
        \KwResult{a subset $K' \subseteq K$ of size at most $(3t_I)^r$ such that $\textbf{s}(K')|_I  = \textbf{s}(K)|_I $}

        Set $U:=K$ and $Q:=\emptyset$, where $Q$ is maintained as an ordered set.

\While{$U\neq \emptyset$}{

        \eIf{$S^+(Q)> t_{I}$}
         {Choose $T\in U$ such that $S^+(Q\cup \{T\})\le  t_{I}+2(r-1)$;\\
          set $Q:=Q\cup \{T\}$ and $U:=U\setminus \{T\}$; \\}
         {let $T\in U$ be an arbitrary copy of $K_r$;\\
          set $Q:=Q\cup \{T\}$ and $U:=U\setminus \{T\}$. \\}


\If{there exists nonempty $Q'\subseteq Q$ such that $\textbf{s}_I(Q') = \textbf{0}$} 
{\label{line:if_zero}
set $U:=U\cup (Q\setminus Q')$ and $Q:=\emptyset$; \\}
         }

        Output $K':=Q$.
\end{algorithm}

In the proof of Lemma~\ref{lem:slack}, we choose the sub-family $K'$ of $K$ by repeated applications of the pigeonhole principle.
Indeed, we track the sum of the absolute values of the entries of the slack vectors restricted to $I$.
More precisely, for any $K_r$-tiling $Q$ on $V$, let $\textbf{s}_I(Q) := \textbf{s}(Q)|_I$. 
Denote by $\textbf{s}_{I,j}(Q)$ the $j$-th coordinate of $\textbf{s}_I(Q)$ with $j\in [|I|]$.
Let $S^+ (Q) := \sum_{j\in [|I|]} |\textbf{s}_{I, j}(Q) |$. Note that for any copy $T'$ of $K_r$ in $K$, we know $S^+(\{T'\}) \leq 2(r-1)$.
We shall let Algorithm~\ref{alg:slack} read the members of $K$ greedily but with a slack constraint in the sense that the partial $K_r$-tiling $Q$ considered by Algorithm~\ref{alg:slack} will never have large $S^+ (Q) $.
When $Q$ is large enough, by the pigeonhole principle we can find a nonempty $K_r$-tiling $Q'$ with $\textbf{s}_I(Q') = \textbf{0}$.
We then remove the members of $Q'$, put the members of $Q\setminus Q'$ back to $U$ and repeat this process.

\begin{proof}[Proof of Lemma~\ref{lem:slack}]
We first prove the case $r\ge 4$.
The lemma is trivial if $(3t_I)^r\ge n/r$. So we may assume that $(3t_I)^r< n/r$.
We execute Algorithm \ref{alg:slack} and add copies of $K_r$ from $U$ to $Q$. If $S^+(Q)\le t_{I}$, then $S^+(Q\cup\{T\})\le t_{I}+2(r-1)$ for any $T\in U$, as $S^+(\{T\}) \leq 2(r-1)$.
First, since each time when we remove a $K_r$-tiling $Q'$ (on lines 10-11), we have $\textbf{s}_I(Q')=\textbf{0}$, at any stage of the algorithm, we have $\textbf{s}_I(Q\cup U)=\textbf{s}_I(K)$, and thus $S^+(Q\cup U)=S^+(K)=t_I$.
Another immediate consequence is that for the final output $K'$, as $U=\emptyset$ and $Q=K'$, we conclude that $\textbf{s}_I(K')=\textbf{s}_I(K)$.

We next show that when $S^+(Q)> t_{I}$ the choice of $T$ as on line 4 of Algorithm \ref{alg:slack} is always possible.
Indeed, suppose there is a moment that $S^+(Q)> t_{I}$ holds but it did not hold in the previous iteration.
Then $t_{I}<S^+(Q)\le t_{I}+2(r-1)$.
Since $S^+(Q\cup U) = t_{I} < S^+ (Q)$, there exists a copy $T$ of $K_r$ in $U$ such that $ S^+ (Q\cup \{T\})< S^+ (Q)\le t_{I}+2(r-1)$, that is, we can always choose $T\in U$ as described in line 4 of Algorithm \ref{alg:slack}. 
Hence throughout the process, we can guarantee $S^+(Q)\le t_{I}+2(r-1)$.


It remains to show that when $Q$ gets large enough, on line 10 we can find the desired $Q'$, which will also upper bound $|K'|$.
Note that each time we execute line 11, we ``restart'' the process by putting the copies of $K_r$ in $Q\setminus Q'$ back to $U$.
Therefore, let $Q$ be a nonempty ordered set of copies of $K_r$ before the execution of line 10, and for $i\in [|Q|]$ let $Q_i$ be the set of the first $i$ copies of $K_r$ in $Q$.
Then for each $i\in [|Q|]$, by line 4, $Q_i$ satisfies that $S^+(Q_i)\in [0, t_I+2(r-1)]$.
Hence 
\begin{itemize}
    \item[($\dagger$)] $\textbf{s}_{I, j}(Q_i)\in [-(t_I+2(r-1)), t_I+2(r-1)]$ for each $i\in [|Q|]$ and $j\in [|I|]$.
\end{itemize}
Therefore, the number of possible values of $\textbf{s}_I(Q_i)$ is at most $(2( t_I+2(r-1))+1)^{|I|}\le (3t_I)^r$.
Moreover, if line 11 is not executed in this iteration of the while loop, then in particular $\textbf{s}_I(Q_i)\neq \textbf{s}_I(Q_j)$ for any $1\le i<j\le |Q|$, we conclude that $|Q|\le (3t_I)^r$.
In particular, we have $|K'|\le (3t_I)^r$.

Finally, let $r=3$ and we shall obtain a refined bound.
First if $r'=2$, then by ($\dagger$), the number of possible values of $\textbf{s}_I(Q_i)$ is at most 
$3t_I$. Indeed, this is trivial if $|I|=1$; for $|I|=2$ it also holds as $\textbf{s}_{I, 1}(Q_i)=-\textbf{s}_{I, 2}(Q_i)$ and there are at most 
$3t_I$ choices for either of them. 

Now suppose $r'=3$.
Again the result is trivial if $|I|=1$.
If $|I|\ge 2$, then note that $t_{[3]}\le 2t_I$, and thus it suffices to show that the 
$K_r$-tiling $K'$ can be chosen to have size at most $4t_{[3]}$. 
Indeed, if $|I|=2$, then by the assumption, we obtain a 
$K_r$-tiling $K'$ of size at most $4t_{[3]}$ such that $\textbf{s}(K')=\textbf{s}(K)$. 
Hence $|K'|\le 4t_{[3]}\le 8t_I$ and $\textbf{s}(K')|_I=\textbf{s}(K)|_I$, that is, $K'$ is a desired $K_r$-tiling.

For $I=[3]$, note that there are three types of index vectors of triangles up to permutation, that is, $(3,0,0)$, $(2,1,0)$ and $(1,1,1)$, whose slack vectors are $(2,-1,-1)$, $(1,0,-1)$ and $(0,0,0)$, respectively.
Suppose $K'$ is minimal subject to the condition $\textbf{s}(K')=\textbf{s}(K)$.
Therefore, $K'$ contains no triangles of index vector $(1,1,1)$, as we may freely remove them without changing $\textbf{s}(K')$.
We now notice a useful fact stated as follows. 
\begin{claim}\label{clm:slack_non_neg}
    If there is $i\in [3]$ such that all slack vectors of triangles in $K'$ have non-negative (non-positive) values on the $i$-th digit, then $|K'|\le 4t_I$.
\end{claim}
\begin{proof}
    Suppose all triangles of $K'$ have slack vectors non-positive in the first digit. Let $K_1'\subseteq K'$ be the triangles in $K'$ whose first digit negative. Then we have $|K_1'|\le t_I$. For $K_2':=K'\setminus K_1'$, note that all triangles in it have the same slack vectors, either all $(0,1,-1)$, or all $(0,-1,1)$.
    Without loss of generality, suppose all of them have slack $(0,1,-1)$, yielding that $\textbf{s}(K_2')=(0,|K_2'|,-|K_2'|)$.
    As $|\textbf{s}_{I, j}(K_1')|\le 2|K_1'|\le 2t_I$ for $j=2,3$, we conclude that $||K_2'|-2t_I|\le t_I$, implying $|K_2'|\le 3t_I$. Therefore, we have $|K'|\le 4t_I$.
\end{proof}

Now we split the discussion according to the types of triangles in $K'$.
Let $\mathcal I $ be the set of the slacks of triangles in $K'$, that is, $\bfi \in \mathcal I$ if and only if $K'$ contains a triangle of slack $\bfi$.

\noindent\textbf{Case 1. $K'$ contains triangles of slack vectors $(2,-1,-1)$ and $(-1,2,-1)$, that is, $(2,-1,-1)$, $(-1,2,-1)\in \mathcal I$.}
We may suppose that $(-1,0,1)$ or $(0,-1,1)$ is in $\mathcal I$, as otherwise we conclude $|K'|\le 4t_I$ by Claim \ref{clm:slack_non_neg}.
Without loss of generality, suppose $(-1,0,1)\in \mathcal I$.
Then by minimality, we know that $K'$ has no triangle with slack $(1,0,-1)$, $(-1,1,0)$ or $(0,-1,1)$ and at most one triangle with slack $(1,-1,0)$.
Therefore, $K'$ contains triangles of slack $(2,-1,-1)$, $(-1,2,-1)$, $(-1,0,1)$ and possibly $(1,-1,0)$ (at most one) and $(0,1,-1)$.
Note that only $(-1,0,1)$ has positive slack on its third coordinate.
Now if $K'$ contains at most $t_I+3$ triangles of slack $(-1,0,1)$, then it contains at most 
$2t_I+3$ triangles of slack being one of $(2,-1,-1)$, $(-1,2,-1)$ and $(0,1,-1)$, yielding 
$|K'|\le 3t_I+7\le 4t_I$.
Otherwise $K'$ contains at least $t_I+4$ triangles of slack $(-1,0,1)$, and thus contribute $-t_I-4$ to the first digit of $\textbf{s}(K')$.
As $K'$ contains at most one triangle of slack $(1,-1,0)$, this implies that $K'$ contains at least two triangles of slack $(2,-1,-1)$.
However, noticing the equation $2(2,-1,-1)+(-1,2,-1)+3(-1,0,1)=(0,0,0)$, we obtain a contradiction with the minimality of $K'$.

\noindent\textbf{Case 2. $(2,-1,-1)\in \mathcal I$ and $(-1,2,-1),(-1,-1,2)\notin \mathcal I$.}
Similar to the previous case, by Claim \ref{clm:slack_non_neg}, we may suppose that $(-1,1,0)\in \mathcal I$ ($\mathcal I$ needs to have $(-1,1,0)$ or $(-1,0,1)$, which are the only slack vectors with negative first coordinate).
Then by minimality of $K'$, we have 
$(1,-1,0), (-1,0,1)\notin \mathcal I$. 
Now, considering the third coordinate, we must have $(0,-1,1)\in \mathcal I$ by Claim \ref{clm:slack_non_neg} and so $(0,1,-1), (1,0,-1)\notin \mathcal I$ by minimality.
In summary, $\mathcal I=\{(2,-1,-1),(-1,1,0),(0,-1,1)\}$ and $K'$ contains at most one triangle of slack $(-1,1,0)$.
Note that $(-1,1,0)$ is the only one with positive coordinate on the second coordinate and $K'$ has at most one such triangle. Then as the other two have negative second coordinate, we conclude that $|K'|\le t_I+2$.

\noindent\textbf{Case 3. $(2,-1,-1), (-1,2,-1),(-1,-1,2)\notin \mathcal I$.}
Clearly in this case $|\mathcal I|\le 3$.
If $|\mathcal I|\le 2$, then by Claim \ref{clm:slack_non_neg}, $|K'|\le 4t_I$.
Otherwise $|\mathcal I|= 3$, and without loss of generality, suppose $\bfv_1,\bfv_2, \bfv_3\in \mathcal I$ (they are all permutations of $(-1,1,0)$).
We indeed have $\bfv_3=\bfv_1+\bfv_2$ and for any non-zero coordinate of $\bfv_3$, the equation $\bfv_3=\bfv_1+\bfv_2$ restricts on this coordinate is either $1=0+1$ or $-1=0+(-1)$.
Therefore, we conclude $|K'|\le 4t_I$ by Claim \ref{clm:slack_non_neg}.

We conclude the proof by symmetry.
\end{proof}

Next we give our color-coding lemma, which slightly extends its original formulation.
Given a partition $(V_1, V_2, \dots, V_{r'})$ of $V(G)$, consider all types of $K_r$, that is, all distinct vertex distributions of $K_r$ on $V(G)$. Recall that there are $y={{r+r'-1}\choose {r'-1}}$ different types of $K_r$. We fix an arbitrary ordering of these $y$ types.
Given a $K_r$-tiling of size $m$, its \emph{pattern} is an $m$-dimensional vector where each digit is in $[y]$ representing its type.
Denote by $\mathcal P_m$ the collection of all patterns of $K_r$-tilings of size $m$.

\begin{lemma}[Color-coding lemma]\label{lem:color coding}
Let $r\geq r'\geq 2$ and $m<n/r$.
Let $G$ be an $n$-vertex graph and $(V_1, V_2, \dots, V_{r'})$ be a partition of $V(G)$.
Given a pattern $P\in \mathcal P_m$.
Then there is an algorithm that decides whether $G$ has a $K_r$-tiling of size $m$ with pattern $P$ in time $O(n^r(\log^2 n)2^{O(rm)})$.
\end{lemma}

\begin{proof}
Let $K_r(G)$ be the set of copies of $K_r$ in $G$ and $K_r^i(G)$ be the set of copies of $i$th-type $K_r$ in $G$, where $i\in [y]$. Let $P=(\kappa_1,\ldots,\kappa_m)$.

The color-coding technique requires an explicit construction of a family of perfect hash functions and we use the one given by Fredman, Kolm\'os, and Szemer\'edi~\cite{FKS}.
For any $n, k\in \mathbb N$ and set $U$ of size $n$,
$\mathcal F$ is a family of functions $f: U\to X$ with $|X|=O(k)$
such that for any $k$-subset $U'$ of $U$
there is a function $f\in \mathcal F$ that is injective on $U'$.
Moreover, the family $\mathcal F$ has size $2^{O(|X|)}\log^2 n$
and any function evaluation on a given input takes $O(1)$ time.
We use such a family $\mathcal F$ with $U=V$
and $X$ being a set of size $O(k)=O(rm)$ where $k=rm$.
We call the elements in $X$ \emph{colors}.

For each function $f\in \mathcal F$,
we will compute families $\mathbf C_1$, $\mathbf C_2$, $\dots$, $\mathbf C_m$,
where for $j \ge 1$ each $\mathcal C \in \mathbf C_j$
is a sequence of $j$ pairwise disjoint $r$-subsets of $X$
that can be realized as the color set via $f$ of a $K_r$-tiling in $G$,
i.e.\ there is an ordered $K_r$-tiling $(H_1,\dots,H_j)$ in $G$
such that $\mathcal C = (f(H_1),\dots, f(H_j)$), where in addition we require that the type of $H_i$ is $\kappa_i$ for $i\in [j]$.
Let $\bigcup \mathcal C := \bigcup_{C \in \,\mathcal C} C$.
We will construct $\mathbf C_j$ such that each $rj$-subset of $X$
that can be realized as the color set of the vertices covered
by a $K_r$-tiling of size $j$ in $G$ occurs exactly once
as $\bigcup \mathcal C$ for some $\mathcal C \in \mathbf C_j$.

We first construct $\mathbf C_1$ by considering each copy $H$ of $K_r$ in $K_r^{\kappa_1}(G)$ sequentially
and adding $f(H)$ to $\mathbf C_1$ if $|f(H)|=r$ and it has not been present yet (recall that $f(H)$ is a set of colors).
Suppose that we have constructed $\mathbf C_j$ for some $j \ge 1$.
Since each $rj$-subset of $X$ is represented at most once, we have $|\mathbf C_j|\le \binom {|X|}{rj}$.
Now we construct $\mathbf C_{j+1}$.
For convenience, let $i:=\kappa_{j+1}\in [y]$.
We consider each
$\mathcal C\in \mathbf C_j$ and $H\in K_r^i(G)$ sequentially,
and add to $\mathbf C_{j+1}$ the sequence
$\mathcal C'$ obtained by appending $f(H)$ to $\mathcal C$,
unless $|f(H)|<r$ or $f(H) \cap \bigcup \mathcal C \ne \emptyset$
or $f(H) \cup \bigcup \mathcal C$ already occurs
as $\bigcup \mathcal C''$ for some $\mathcal C'' \in \mathbf C_{j+1}$.
This process continues until $\mathbf C_{m}$ has been constructed.
Hence the running time of this process is
\[
|\mathcal F|\cdot O(m\cdot 2^{2|X|}\cdot |K_r(G)|)
= 2^{O(rm)}\log^2 n \cdot O(m\cdot 2^{O(rm)}\cdot n^r) = O(n^r (\log^2 n) 2^{O(rm)}).
\]
If the above process returns $\mathbf C_{m}=\emptyset$ for each $f \in \mathcal F$,
then
$G$ does not contain a $K_r$-tiling of size $m$ by the defined property of the perfect hash family $\mathcal F$.
Otherwise, suppose we obtain $\mathcal C^* \in \mathbf C_{m}$ for some $f \in \mathcal F$.
For each $j\in [m]$, let $\{x_{(j-1)r+1}, \dots, x_{jr}\}$ be the $j$-th $r$-set in $\mathcal C^*$.
Based on our construction, $G[f^{-1}(x_{(j-1)r+1}), \dots, f^{-1}(x_{jr})]$
induces at least one copy of $K_r$, denoted by $H_j^{\kappa_j}$ for each $j$, which can be found in time $O(|K_r(G)|)=O(n^r)$.
Then $H_1^{\kappa_1},\dots,H_m^{\kappa_m}$ form a $K_r$-tiling of size $m$ with pattern $P$ in $G$.
The total running time is $O(n^r (\log^2 n) 2^{O(rm)})$.
\end{proof}

\section{Extremal case}\label{sec: ext case}
In this section we study the extremal case.
Throughout the proof, it suffices to deal with
\[
c = \begin{cases}
    o(n^{1/r^{r+3}}), \text{ if } r\ge 4, \\
    o(n), \hfill\text{ if } r=3,
\end{cases} 
\]
as otherwise, 
Theorem~\ref{thm:fellows}  implies that we can determine whether $G$ contains a $K_r$-factor in time $O(n^r+2^{O(n)})=2^{O(n)}$, which is $2^{O(c)}$ for $r=3$ and $2^{O(c^{r^{r+3}})}$ for $r\ge 4$, as desired.
Moreover, throughout this section, we work with the following parameters satisfying the hierarchy
\[
1/n\ll \gamma \ll \mu \ll \lambda \ll 1/r,
\]
which we will omit in the subsequent statements of lemmas.

In the proof of the non-extremal case, when the proof fails, we detect a sparse $\frac nr$-set.
However, in the proof of the extremal case, to get a clear picture, we need to know how many vertex-disjoint such sparse sets $G$ contains.
For this goal our first step is to identify these sparse sets, using the regularity method, which is known to be efficient for hunting rough structures of the graph.

\subsection{Find sparse sets}
The following lemma is the main result of this section.
\begin{lemma}\label{lem:find sparse sets}
Let $G$ be an $n$-vertex graph with $\delta(G)\ge (1-1/r)n-c$.
If $G$ has at least one $\gamma$-sparse $\frac{n}{r}$-set, then by Algorithm \ref{alg:SparseSet} we can find $s$ pairwise disjoint $\mu$-sparse $\frac{n}{r}$-sets $A_1,A_2,\ldots,A_s$ of $G$ in time $O(n^{2.376})$ such that $B:=V(G)\setminus \cup_{i=1}^sA_i$ has no $\lambda$-sparse $\frac{n}{r}$-set of $G$, where $1\le s \le r$.
\end{lemma}

\begin{algorithm}
\caption{Find sparse sets}\label{alg:SparseSet}
\KwData{an $n$-vertex graph $G$ with $\delta(G)\ge (1-1/r)n-c$ and a $\gamma$-sparse $\frac{n}{r}$-set $A_1$ in $G$}
\KwResult{sparse $\frac{n}{r}$-sets of $G$}
Choose constants $1/n\ll\g \ll \g_2\ll \e_2,d_2\ll\g_3\ll \cdots \ll \g_{r}\ll \e_{r},d_{r}\ll \g_{r+1}\ll 1/r$;

\For{$i=2$ to $r$}{

Let $G_i:=G\setminus \cup_{j=1}^{i-1}A_j$ and $n_i:=|G_i|=(r-i+1)\frac{n}{r}$;

Apply Lemma \ref{lem:reg} on $G_i$ with parameters $\e_i$ and $d_i$ to obtain a partition $V(G_i)=V^{i}_0\cup V^{i}_1\cup\ldots \cup V^{i}_{k_{i}}$ for some $1/\e_i \le k_i\le M(\e_i)$;

Let $R_{i}$ be the reduced graph for this partition and so $|R_i|=k_i$;

\eIf{$R_i$ has no $\left(\frac{8\g_{i}}{\e_i^2d_i}+\frac{3\e_i}{r-i+1}\right)$-sparse $\frac{k_i}{r-i+1}$-set}
         {$G_i$ has no $\gamma_{i}$-sparse $\frac{n_i}{r-i+1}$-set by Fact \ref{lem:sparse set2};\\
          Output $A_1,A_2,\ldots, A_{i-1}$ and halt. \\
         }
         {$G_i$ has $\left(\frac{8\g_{i}}{\e_i^2d_i}+\frac{3\e_i}{r-i+1}+2\e_i+d_i\right)$-sparse $\frac{n_i}{r-i+1}$-set, denoted by $A_{i}$, by Fact \ref{lem:sparse set1}. \\
         }

}
\end{algorithm}

\begin{proof}
We choose the constants as in Algorithm~\ref{alg:SparseSet} and let $\mu_1=\g$. For $i\in [2,r]$, let $\mu_i:=\frac{8\g_{i}}{\e_i^2d_i} + \frac{3\e_i}{r-i+1} + 2\e_i + d_i$.
Then we have $\g_i\ll \mu_i\ll \g_{i+1}$ as $\g_i\ll \e_i,d_i\ll \g_{i+1}$.
Let $G$ be a graph as assumed in the lemma and we execute Algorithm~\ref{alg:SparseSet} to $G$.
Note that if $A_{i}$ is defined for any $i\in [2,r]$, then $A_{i}$ is a $\mu_i$-sparse $\frac{n}{r}$-set, as $n_i=(r-i+1)\frac{n}{r}$ and $|A_i|=\frac{n_i}{r-i+1}=\frac{n}{r}$.

Let $s$ be the number of sparse sets that Algorithm \ref{alg:SparseSet} outputs. Let $\mu=\mu_s$ and $\lambda=\g_{s+1}$.
Suppose that Algorithm \ref{alg:SparseSet} halts when $i\le r-1$. Then $s=i-1$ and we can find $s$ pairwise  disjoint $\mu$-sparse $\frac{n}{r}$-sets $A_1,A_2,\ldots,A_s$ such that $B:=G\setminus \cup_{j=1}^{s}A_j$ has no $\lambda$-sparse $\frac{n}{r}$-set.

Otherwise suppose that Algorithm \ref{alg:SparseSet} halts with $i=r$. 
Then either $s=r-1$ or $s=r$. 
If $s=r-1$, then we can find $r-1$ pairwise  disjoint $\mu$-sparse $\frac{n}{r}$-sets $A_1,A_2,\ldots,A_{r-1}$ such that $B:=G\setminus \cup_{j=1}^{r-1}A_j$ is not a $\lambda$-sparse $\frac{n}{r}$-set. If $s=r$, then we can find $r$ pairwise disjoint $\mu$-sparse $\frac{n}{r}$-sets $A_1,A_2,\ldots,A_r$ of $G$.

By Lemma \ref{lem:reg} and noting that all $k_i\le M(\eps_i)$ are constants, the running time of Algorithm \ref{alg:SparseSet} can be bounded by $\sum_{i=2}^r\left(O(n^{2.376})+{k_i\choose {k_i/(r-i+1)}}\right)=O(n^{2.376})$.
\end{proof}


We need the following result when $s=r-2$ and $B$ is close to the union of two disjoint cliques.

\begin{lemma}\label{lem:odd comp}
Let $B'\subseteq V(G)$ be a subset such that $|B'|\ge 2n/r-\mu^{1/3} n$ and $|B'|$ is even.
If $G[B']$ has no perfect matching, then $G[B']$ is disconnected with two odd components both of which are almost complete.
\end{lemma}

\begin{proof}
Note that $\delta(G[B'])\ge |B'| - (n/r+c)\ge |B'| - (|B'|/2+\mu^{1/3} n) \ge (1/2-\mu^{1/3}r)|B'|$.
Since $G[B']$ has no perfect matching, by Tutte's theorem, there exists $S\subseteq B'$ such that the subgraph $G[B'\setminus S]$ has $x> |S|$ odd components.
Suppose that $|S|\neq 0$. Then the smallest odd component has at most $(|B'|-|S|)/|S|$ vertices, and by the minimum degree condition, we have $(|B'|-|S|)/|S|+|S| > (1/2-2\mu^{1/3}r)|B'|$, and thus either $|S| > (1/2-3\mu^{1/3}r)|B'|$ or $|S|<1/(1/2-3\mu^{1/3}r)$. 
The former implies that $G[B'\setminus S]$ contains an independent set of size $x> |S|> (1/r-7\mu^{1/3})n$ (by taking one vertex from each odd component), which implies a $\lambda$-sparse $\frac{n}{r}$-set in $G$, a contradiction.
The latter says that $|S| \le 2$ and by the minimum degree condition, $G[B'\setminus S]$ has at most two components.
If $|S|=2$, then $G[B'\setminus S]$ has $x> 2$ odd components, a contradiction.
If $|S|=1$, then $G[B'\setminus S]$ has two odd components, which contradicts the parity of $|B'|$.
Hence, we have $|S|=0$ and $G[B']$ has at least one odd component.
By the minimum degree condition, $G[B']$ is disconnected with two odd components both of which have order at least $(1/2-\mu^{1/3}r)|B'|$ and are almost complete.
\end{proof}



Recall that by Lemma \ref{lem:find sparse sets}, we have found $s$ pairwise disjoint $\mu$-sparse $\frac{n}{r}$-sets $A_1,A_2,\ldots,A_s$ of $G$ such that $B:=V(G)\setminus \cup_{i=1}^sA_i$ has no $\lambda$-sparse $\frac{n}{r}$-set of $G$, where $1/n\ll\mu\ll \lambda\ll 1/r$.
In the rest of the proof, there is a tiny difference between the cases $1\le s\le r-1$ and $s=r$.
We will first prove the former case with full details in Sections \ref{sec: ideal case}-\ref{sec: obtain K_r-s-factor}, and explain how to modify the approach to accommodate the latter in Section \ref{sec: s=r}. 

\subsection{Ideal case}\label{sec: ideal case}

We first present our ideal case which allows us to find a $K_r$-factor easily. Hence in Sections \ref{sec: balancedness}-\ref{sec: obtain K_r-s-factor} we would like to get rid of the atypical vertices while keeping the correct ratio to achieve this ideal case. 
Note that in the following lemma we do not require that for example $A_i'\subseteq A_i$ or $B'\subseteq B$ -- indeed, this will not be the case as we need to adjust these sets slightly.

\begin{lemma}\label{lem: ideal}
    Suppose $s<r$ and we have $A_1', \dots, A_s', B'$ such that 
    \begin{itemize}
        \item $|A_1'|=\cdots=|A_s'|=|B'|/(r-s)\ge n/r-\mu^{1/5}n$;
        \item for every vertex $v\in B'$ and any $i\in [s]$, $d_{A_i'}(v)\ge |A_i'|-\mu^{1/5}n$;
        \item for any $i\in [s]$ and every vertex $u\in A_i'$, $d_{A_j'}(u)\ge |A_j'|-\mu^{1/5}n$ for any $j\in [s]\setminus \{i\}$ and $d_{B'}(u)\ge |B'|-\mu^{1/5}n$;
        \item $G[B']$ has a $K_{r-s}$-factor.
    \end{itemize}
Then in time $O(sn^4)$ we can find a $K_r$-factor of $G[A_1'\cup \dots\cup A_s'\cup B']$. 
\end{lemma}
\begin{proof}
Note that a $K_{r-s}$-factor of $G[B']$ has size $\frac{|B'|}{r-s}$.
We arbitrarily fix a $K_{r-s}$-factor $\mathcal{T}$ of $G[B']$.
We contract all members of $\mathcal{T}$ individually to a set of disjoint vertices, denoted by $U$. 
For each $u\in U$, let $K_{r-s}^u$ be the copy of $K_{r-s}$ corresponding to $u$.
Next we construct an auxiliary bipartite graph $H_s$ with two parts $U$ and $A_s'$. 
Note that $|U| = |A_s'| = \frac{|B'|}{r-s} $.
For any $u\in U$ and $v\in A_s'$, $uv\in E(H_s)$ if and only if $v$ is a common neighbor of vertices in $K_{r-s}^u$. Hence $d_{H_s}(u)\ge |A_s'|-(r-s)\mu^{1/5}n > |A'_s|/2$ and $d_{H_s}(v)\ge \frac{|B'|}{r-s}-\mu^{1/5}n > |A'_s|/2$, as every vertex in $K_{r-s}^u$ has at most $\mu^{1/5}n$ non-neighbors in $A_s'$ and $v$ has at most $\mu^{1/5}n$ non-neighbors in $B'$.
Thus $\delta(H_s)>|A'_s|/2$ and so there is a perfect matching in $H_s$, i.e. a matching of size $|A_s'|$. We can find the matching in time $O(n^4)$ by the Edmonds algorithm, which corresponds to a $K_{r-s+1}$-factor of size $|A_s'|$ in $G[A_s'\cup B']$. (Note that if $r-s=1$, then $H_s$ is $G[A_s'\cup B']$.)

Iterating this procedure $s$ times we obtain a $K_{r}$-factor of size $|A_1'|$ in $G[A_1'\cup \dots\cup  A_s'\cup B']$.
\end{proof}

For the rest of the proof, our goal is to remove a small $K_r$-tiling from $G$ so that we can apply Lemma~\ref{lem: ideal}, and on the other hand, when $G$ has no $K_r$-factor, we need to show that such a (small) $K_r$-tiling does not exist.

\subsection{Adjusting the partition and more definitions}\label{sec: balancedness}
Note that the current partition $A_1,\dots, A_s$, $B$ only captures the rough structure of $G$.
Motivated by the ideal case, we need to further classify the vertices of $G$.
In this section, we move a small number of vertices and obtain a new partition $A_1^1,A_2^1,\ldots,A_s^1, B^1$ such that vertices in each $A_i^1$ have small degrees inside $A_i^1$ and have large degrees outside 
$A_i^1$. 
However, $A_1^1,A_2^1,\ldots,A_s^1,B^1$ may be unbalanced. Then by removing a $K_r$-tiling of small size from $G$, they will become balanced again (Lemmas~\ref{lem:balancedness} and \ref{lem:balancedness2}).

The following definitions are motivated by Lemma~\ref{lem: ideal}. 

\begin{definition}\label{def: Ai-good vtx}
For any $i\in [s]$, a vertex $v\in V(G)$ is called an \emph{$A_i$-good vertex} if $d_{A_i}(v)< 2\mu^{1/3}n$.
\end{definition}

\begin{fact}\label{fact: properties of Ai-good vtx}
For each $i\in[s]$, we have the following properties:
 \begin{enumerate}
   \item [(1)] the number of not $A_i$-good vertices in $A_i$ is less than $\mu^{2/3}n$;
   \item [(2)] the number of $A_i$-good vertices in $V(G)\setminus A_i$ is less than $4\mu r n+2c$;
   \item [(3)] any vertex $v\in V(G)$ cannot be both $A_i$-good and $A_j$-good for $i\neq j$.
 \end{enumerate}
\end{fact}
\begin{proof}
(1) Suppose that the number of not $A_i$-good vertices in $A_i$ is at least $\mu^{2/3}n$. Since $d_{A_i}(v)\ge 2\mu^{1/3}n$ for any such vertex $v$, we have $e(G[A_i])\ge \mu n^2$, a contradiction.

(2) Since $e(G[A_i,V(G)\setminus A_i])\ge (n-n/r-c)\cdot|A_i|-2e(G[A_i])>|V(G)\setminus A_i||A_i|-2\mu n^2-|A_i|c$, there are at most $2\mu n^2+|A_i|c$ non-edges between $A_i$ and $V(G)\setminus A_i$. Moreover, an $A_i$-good vertex in $V(G)\setminus A_i$ contributes at least $|A_i|-2\mu^{1/3}n$ non-edges between $A_i$ and $V(G)\setminus A_i$, so the number of $A_i$-good vertices in $V(G)\setminus A_i$ is at most $\frac{2\mu n^2+|A_i|c}{|A_i|-2\mu^{1/3}n}=\frac{2\mu rn+c}{1-2\mu^{1/3}r}<4\mu r n+2c$.

(3) Since $\delta(G)\ge (1-1/r)n-c$, any vertex $v\in V(G)$ cannot be both $A_i$-good and $A_j$-good for $i\neq j$.
\end{proof}

Now we define the vertex partition $A_i^1:=\{v: v~\text{is an}~A_i\text{-good vertex}\}$ for each $i\in [s]$ and $B^1:=V(G)\setminus \cup_{i=1}^{s}A_i^1$.
Then $A_1^1,A_2^1,\ldots,A_s^1,B^1$ may be unbalanced.
We assume that $|A_i^1|=n/r+t_i$ for each $i\in [s]$ and $|B^1|=(r-s)n/r+t_{s+1}$.
Note that $\sum_{i=1}^{s+1}t_i=0$.
By Definition \ref{def: Ai-good vtx} and Fact \ref{fact: properties of Ai-good vtx}, we have the following result.

\begin{lemma}\label{lem: move vtx}
The sets $A_i^1$ for each $i\in [s]$ and $B^1$ have the following properties:
\begin{enumerate}
   \item [(1)] $\delta(G[A_i^1])\ge t_i-c$ and $\delta(G[B^1])\ge (r-s-1)n/r+t_{s+1}-c$;
   \item [(2)] $-\mu^{2/3}n<t_i<4\mu rn+2c$ for each $i\in [s]$ and $-4s\mu rn-2sc<t_{s+1}<s\mu^{2/3}n$;
   \item [(3)] $\Delta(G[A_i^1])< 2\mu^{1/3}n+4\mu rn+2c$;
   \item [(4)] every vertex $v\in B^1$ satisfies $d_{A_i}(v)\ge 2\mu^{1/3}n$ for each $i\in [s]$;
   \item [(5)] $B^1$ has no $\frac{\lambda}{2}$-sparse $\frac{n}{r}$-set of $G$.
 \end{enumerate}
\end{lemma}

\begin{proof}
(1) Since $\delta(G)\ge (1-1/r)n-c$, for any $v\in V(G)$, it has at most $\frac n r +c$ non-neighbors.
Thus $\delta(G[A_i^1])\ge |A_i^1| - (\frac n r +c) \ge t_i-c$ and $\delta(G[B^1])\ge |B^1| - (\frac n r +c)\ge (r-s-1)n/r+t_{s+1}-c$.

(2) By Fact~\ref{fact: properties of Ai-good vtx} (1) and (2), for each $i\in[s]$, we have $-\mu^{2/3}n<t_i<4\mu rn+2c$. Then \[-4s\mu rn-2sc<t_{s+1} = -\sum_{i=1}^s t_i<s\mu^{2/3}n.\]

(3) By Definition~\ref{def: Ai-good vtx} and Fact~\ref{fact: properties of Ai-good vtx} (2), we have
\[\Delta(G[A_i^1]) < 2\mu^{1/3}n+4\mu rn+2c.\]

(4) By Definition~\ref{def: Ai-good vtx} and the definition of $B^1$, any vertex $v\in B^1$ satisfies $d_{A_i}(v)\ge 2\mu^{1/3}n$ for each $i\in [s]$.

(5) By Lemma~\ref{lem:find sparse sets}, we know that $B$ has no $\lambda$-sparse $\frac nr $-set of $G$, that is, for any $W\subseteq B$ with $|W|=\frac nr$, $e(G[W])\ge \lambda n^2$. 
Then for any $\frac nr$-set $U\subseteq B^1$, we have
\[e(G[U])\ge \lambda n^2 - 2 |U| \cdot (s\mu^{2/3}n)\ge \lambda n^2 - 2  s\mu^{2/3}n^2/r \ge \frac \lambda 2 n^2 ,\]
as $\mu \ll \lambda$.
Thus $B^1$ has no $\frac{\lambda}{2}$-sparse $\frac{n}{r}$-set of $G$.
\end{proof}

By Lemma \ref{lem: move vtx} (4), every vertex $v\in B^1$ satisfies $d_{A_i}(v)\ge 2\mu^{1/3}n$ for each $i\in [s]$. 
However, some vertices in $B^1$ may have many non-neighbors in some $A_i$, that is, they violate the property required in Lemma~\ref{lem: ideal}. 
Hence we mark such vertices as bad and try to cover them by a $K_r$-tiling (Lemma~\ref{lem:clean bad vtx-pre}) so that the remaining part of $G$ can be in the ideal case. Lemma~\ref{lem:clean bad vtx-pre} can cover only a small number of such bad vertices. Fact~\ref{fact: properties of B1-bad vtx} shows that the number of such bad vertices is indeed small.



\begin{definition}\label{def: B1-bad vtx}
A vertex $v\in B^1$ is called \emph{bad} if $d_{A_i}(v)< |A_i|-2\mu^{2/5}n$ for some $i \in [s]$.
\end{definition}

\begin{fact}\label{fact: properties of B1-bad vtx}
There are at most $s(2\mu^{3/5}+\mu^{2/3})n$ bad vertices in $B^1$.
\end{fact}
\begin{proof}
First we consider the bad vertices in $B\cap B^1$. 
For each $i\in[s]$, note that $e(G[A_i,B])\ge |A_i|(\delta(G) - (s-1)n/r) - 2e(G[A_i])\ge |A_i|((r-1)n/r-c - (s-1)n/r) - 2\mu n^2$ and $|B|=(r-s)n/r$.
So we have $e(G[A_i,B])\ge |A_i||B|-2\mu n^2-|A_i|c$, 
and thus there are at most $2\mu n^2+|A_i|c$ non-edges between $A_i$ and $B$. 
For a bad vertex $v$ in $B$, since $d_{A_i}(v)< |A_i|-2\mu^{2/5}n$, it contributes at least $2\mu^{2/5}n$ non-edges between $A_i$ and $B$. So the number of bad vertices in $B$ is at most $s\cdot\frac{2\mu n^2+|A_i|c}{2\mu^{2/5}n}<s\mu^{3/5}n+\frac{s\mu^{3/5}n}{2r}$ as $c\ll \mu n$.
Moreover, as $|B^1\setminus B|< s\mu^{2/3}n$ by Fact~\ref{fact: properties of Ai-good vtx} (1), there are at most $s\mu^{2/3}n$ bad vertices in $B^1\setminus B$. 
Therefore, there are at most $s\mu^{3/5}n+\frac{s\mu^{3/5}n}{2r}+s\mu^{2/3}n\le s(2\mu^{3/5}+\mu^{2/3})n$ bad vertices in $B^1$.
\end{proof}

Let $B^1_b$ be the set of bad vertices in $B^1$ and $B^1_g:=B^1\setminus B^1_b$.
Then $|B^1_b|\le 3s\mu^{3/5}n$ and $|B^1_g|\ge (r-s)\frac{n}{r}+t_{s+1}-3s\mu^{3/5}n \ge (r-s)\frac{n}{r} - 4s\mu^{3/5}n$. 
The next lemma gives a property of $B^1_g$ that will be used in the proof of Lemma~\ref{lem:balancedness2}.

\begin{lemma}\label{lem:Kr-s+1}
If $s<r$, then there are at least $\mu^{7/12}n^{r-s+1}$ copies of $K_{r-s+1}$ in $G[B^1_g]$.
\end{lemma}
\begin{proof}
The minimum degree condition of $G$ says that every vertex has at most $n/r+c$ non-neighbors, which is also true for $G[B^1_g]$. 
By $|B_g^1|\ge (r-s)\frac{n}{r} - 4s\mu^{3/5}n$, we get
\begin{align*}
\frac{\frac nr+c}{|B_g^1|} \le \frac{\frac nr+\mu n}{(r-s)\frac{n}{r} - 4s\mu^{3/5}n} = \frac{1+r\mu}{r-s - 4rs\mu^{3/5}} < \frac1{r-s} + 5rs \mu^{3/5}
\end{align*}
Thus, we have $\delta(G[B^1_g])\ge (1-\frac{1}{r-s}-5rs\mu^{3/5})|B^1_g|$.
Moreover, $B^1_g$ has no $\frac{\gamma_{s+1}}{4}$-sparse $\frac{|B^1_g|}{r-s}$-set, as $\mu\ll \g_{s+1}$.
Then we apply Fact \ref{fact: counting Kr in graph with no sparse set} on $G[B^1_g]$ with $k=r-s+1$, $\varepsilon=5rs\mu^{3/5}$ and $\g_H=\frac{\gamma_{s+1}}{4}$. Since $\mu\ll \g_{s+1}$ and $c = o(n^{1/r^{r+3}})$, there are at least $\mu^{7/12}n^{r-s+1}$ copies of $K_{r-s+1}$ in $G[B^1_g]$.
\end{proof}

\subsection{Admissibility and the characterization lemma}\label{subsec:make balanced} 
The following definition is key to our proofs.
Recall that $|A_i^1|=n/r + t_i$ for each $i\in [s]$,  $|B^1| = (r-s)n/r + t_{s+1}$ and $\sum_{i=1}^{s+1} t_i =0$.
\begin{definition}\label{def:admissible}
Consider the partition $(A_1^1,\dots,A_s^1, B^1)$ and a set $I\subseteq [s+1]$. 
\begin{enumerate}
    \item A $K_r$-tiling $F$ is called \emph{$I$-admissible} if $\textbf{s}(F)|_{i}=t_{i}$ for each $i\in I$ and 
    $|F|\le g(t_I)$.
    \item A $K_r$-tiling $F$ is called \emph{strong admissible} if $\textbf{s}(F)|_{i}=t_{i}$ for each $i\in [s+1]$, 
    $|F|\le g(t_{[s+1]})+2$ and $F$ contains a copy of $K_r$ which contains an odd number of vertices of $B^1$.
\end{enumerate}
 
\end{definition}

The last thing we need to introduce is a set $I^*$ of indices for small slacks.
To distinguish ``large" and ``small", we use Algorithm \ref{alg:t_i}.
It divides $r'$ real numbers into two groups, those with large absolute values and those with small absolute values, in a recursive manner.
Note that the running time of Algorithm \ref{alg:t_i} is $O(r')$.

\begin{algorithm}
\caption{Find small $|t_i|$'s}\label{alg:t_i}
\KwData{$r'$ real numbers, where $2\le r'\le r$, $t_1\ge\cdots\ge t_{i_1}>\frac{c}{0.99}\ge t_{i_1+1}\ge \cdots \ge t_{i_2}\ge -\frac{c}{0.99} > t_{i_2+1}\ge \cdots\ge t_{r'}$}
\KwResult{an index set $I\subseteq [r']$}

Let $j:=i_1$, $\ell:=i_2+1$ and $I:=[i_1+1,i_2]$.

\While{$I \neq [r']$}{

\eIf{$t_j\le \max\{g(rt_{j+1}), |g(rt_{\ell-1})|\}$, or $|t_{\ell}|\le \max\{g(rt_{j+1}), |g(rt_{\ell-1})|\}$}
         {\textbf{Case 1}: if $t_j\le \max\{g(rt_{j+1}), |g(rt_{\ell-1})|\}$ and $|t_{\ell}|\le \max\{g(rt_{j+1}), |g(rt_{\ell-1})|\}$, \\then $I:=I\cup \{j,\ell\}$, $j:=j-1$ and $\ell:=\ell+1$;\\
          \textbf{Case 2}: if $t_j\le \max\{g(rt_{j+1}), |g(rt_{\ell-1})|\}$ but $|t_{\ell}|> \max\{g(rt_{j+1}), |g(rt_{\ell-1})|\}$, \\then $I:=I\cup \{j\}$ and $j:=j-1$;\\
          \textbf{Case 3}: if $|t_{\ell}|\le \max\{g(rt_{j+1}), |g(rt_{\ell-1})|\}$ but $t_j> \max\{g(rt_{j+1}), |g(rt_{\ell-1})|\}$, \\then $I:=I\cup \{\ell\}$ and $\ell:=\ell+1$;\\
         }
         {Output $I$ and halt. \\
         }
}
\end{algorithm}

Before applying Algorithm~\ref{alg:t_i}, we sort $t_1,t_2,\ldots,t_{s+1},\frac{c}{0.99},-\frac{c}{0.99}$ by size and obtain that $t_{\tau(1)}\ge\cdots\ge t_{\tau(i_1)}>\frac{c}{0.99}\ge t_{\tau(i_1+1)}\ge \cdots \ge t_{\tau(i_2)}\ge -\frac{c}{0.99} > t_{\tau(i_2+1)}\ge \cdots\ge t_{\tau(s+1)}$ in time $O(r^2)$.
Let $x_i:=t_{\tau(i)}$ for each $i\in [s+1]$.
Then we apply Algorithm~\ref{alg:t_i} with $r'=s+1$ on $x_1,\ldots,x_{s+1}$ and output $I=[i_1',i_2']$ in time $O(s)$, with the following properties.
Let  $ I^*:=\{\tau(i):i\in I\}$, and so $t_{I^*} = \sum_{i\in I^*} |t_i|$.
We have
\begin{enumerate}
[label=(I\arabic*)]
\item for $i\in I^*$, we have $|t_i|\le(3 rc/0.99)^{r^{r+1}}=(rc/0.33)^{r^{r+1}}$ for $r\ge 4$ and $|t_i|\le (8r)^rc/0.99$ for $r=3$,
\label{item:I1}
\item for $j\in [r']\setminus I^*$, $|t_j|\ge \max\{g(rt_{\tau(i_1')}), |g(rt_{\tau(i_2')})|\}$ and so $|t_j|\ge\max\{g(t_{I^*}),c/0.99\}$, \label{item:I2}
\item $t_{I^*}\le |I|\cdot(rc/0.33)^{r^{r+1}}\le r(rc/0.33)^{r^{r+1}}$ for $r\ge 4$ and $t_{I^*}\le |I|\cdot (8r)^rc/0.99\le (8r)^rrc/0.99$ for $r=3$, and \label{item:I3}
\item if $|I|\ge s$, then $I=I^*=[s+1]$.\label{item:I4}
\end{enumerate}
The last property can be seen easily: since all $t_i$ sum to 0, if at the execution of the algorithm we have $I=[s]$ or $I=[2,s+1]$, then $I$ will also swallow the last index because $|t_{\tau(1)}|\le s |t_{\tau(s+1)}|$ and $|t_{\tau(s+1)}|\le s |t_{\tau(1)}|$.

Now we are ready to state our characterization result.

\begin{lemma}
[Characterization for $s<r$]
\label{lem:charact}
Suppose $\mu>0$ and $1/n\ll \mu$.
Given an $n$-vertex graph $G$ with $\delta(G)\ge (1-1/r)n-c$, let $(A_1^1,\dots,A_s^1,B^1)$ be the partition defined in Section \ref{sec: balancedness} and suppose $s<r$.
Let $I^*\subseteq [s+1]$ and $B_b^1\subseteq B^1$ be defined as above.
Then the following holds.
\begin{enumerate}
\item If $r-s\neq 2$, $I^*\neq [s+1]$, $B_b^1\neq \emptyset$, or $B^1$ is not disconnected with two odd components, then $G$ has a $K_r$-factor if and only if $G$ contains an $I^*$-admissible $K_r$-tiling.
\item Otherwise, $G$ has a $K_r$-factor if and only if $G$ contains a strong admissible $K_r$-tiling.
\end{enumerate}
Moreover, the existence of such $K_r$-tiling in $G$ can be checked by an algorithm that runs in time $n^{O(1)}2^{O(c^{{r}^{r+2}})}$ for $r\ge 4$ (resp. $n^{O(1)}2^{O(c)}$ for $r=3$), and if $G$ has a $K_r$-factor, it also outputs a $K_r$-factor.
\end{lemma}

Lemma~\ref{lem:charact} details the hidden information of Lemma~\ref{lem:char} for the case $0<s<r$ (see more discussions in the conclusion section).

The next two subsections are devoted to the proof of Lemma~\ref{lem:charact}.


\subsection{Auxiliary lemmas}\label{sec: Auxiliary lemmas}
We first present and prove three auxiliary results, Lemmas~\ref{lem:balancedness}--\ref{lem:clean bad vtx-pre}.
The first two are used to make the partition $A_1^1,A_2^1,\ldots,A_s^1$, $B^1$ balanced and the last one is for putting bad vertices into a small $K_r$-tiling.
We remark that parts with small slacks (whose indices are in $I^*$) and large slacks are quite different, so we separate them into two lemmas.

\begin{lemma}\label{lem:balancedness}
There is an algorithm with running time $n^{O(1)}2^{O(c^{{r}^{r+2}})}$ for $r\ge 4$ (resp. $n^{O(1)}2^{O(c)}$ for $r=3$) which either outputs an $I^*$-admissible $K_r$-tiling $K^1$
or returns $(A_1^1,A_2^1,\ldots,A_s^1,B^1)$ as a certificate showing that $G$ has no $I^*$-admissible $K_r$-tiling (and thus also no $K_r$-factor). 
\end{lemma}

\begin{lemma}\label{lem:balancedness2}
Suppose that $I^*\neq [s+1]$ and Lemma \ref{lem:balancedness} outputs an $I^*$-admissible $K^1$. Then we can always find an $[s+1]$-admissible $K_r$-tiling $K$ of size at most $r^2\mu^{2/3}n$ in time $O(n^r)$ such that $K^1\subseteq K$ and $V(G)\setminus V(K)$ has a partition $A_1^2,\ldots,A_s^2,B^2$ with $|A_1^2|=\cdots=|A_s^2|=\frac{|B^2|}{r-s}\ge n/r-r^2\mu^{2/3}n$, $A_i^2\subseteq A_i^1$ for each $i\in [s]$ and $B^2\subseteq B^1$. 
In particular, if $|B^1\setminus V(K^1)|-(r-s)(n/r-|K^1|) > 0$, then $K$ contains a copy of $K_r$ that has exactly $r - s + 1$ vertices of $B^1$.
\end{lemma}

\begin{lemma}\label{lem:clean bad vtx-pre}
Given any bad vertex $v\in B^1$ and any subset $W \subseteq  V(G)$ with $|W|\le r^3\mu^{3/5}n$, we can find a balanced copy $H$ of $K_r$ containing $v$ such that $V(H)\subseteq V(G)\setminus W$ in time $O(n)$.
\end{lemma}



We first prove Lemma \ref{lem:balancedness}.
The key property is~\ref{item:I3}, which allows us to use the color-coding method (Lemma \ref{lem:color coding}) to decide whether there is an $I^*$-admissible $K_r$-tiling in the desired running time.
This is to say, such a search may not be possible for $[s+1]$-admissible $K_r$-tilings when there are very large slacks (could be linear in $n$), however, we shall see later that for large slacks a direct and easy argument exists.
On the other hand, by Lemma \ref{lem:slack}, if such a $K_r$-tiling does not exist, then $G$ has no $K_r$-factor. 

\begin{proof}[{Proof of Lemma \ref{lem:balancedness}.}]
We search for $K_r$-tilings $F$ of size at most $g(t_{I^*})$ in $G$ such that $\textbf{s}(F)|_{i}=t_{i}$ for each $i\in I^*$.
By~\ref{item:I3}, \[g(t_{I^*})=g(t_{I^*})\le (3r(rc/0.33)^{r^{r+1}})^r\le ((r^2c/0.11)^{r^{r+1}})^r=(r^2c/0.11)^{r^{r+2}}\]
for $r\ge 4$, and \[g(t_{I^*})=8t_{I^*}\le (8r)^{r+1}c/0.99\] for $r=3$. (As mentioned at the beginning of Section \ref{sec: ext case}, we only consider $c=o(n^{1/{r^{r+3}}})$ for $r\ge 4$ and $c=o(n)$ for $r=3$. Hence we have $g(t_{I^*})=o(n)$ for $r\ge 3$.)
By Lemma \ref{lem:slack}, if there is no such $K_r$-tiling, then there is no $K_r$-factor in $G$. To find such a $K_r$-tiling or a certificate showing that none exists, we do the following two steps:

\begin{itemize}
  \item \textbf{Step 1.} Search $K_r$-tilings of size at most $g(t_{I^*})$ in $G$.

  \item \textbf{Step 2.} If Step 1 returns some $K_r$-tilings, then we check their slacks.
\end{itemize}

\noindent \textbf{Step 1.} 
Recall from Definition~\ref{def:index_vector} that if we consider the partition $(A_1^1,\dots, A_s^1, B^1)$, then there are $y:={r+s\choose s}$ types of copies of $K_r$ in $G$. 
We fix an arbitrary ordering of these $y$ types. 
Let $x_i$ be the number of copies of $i$th-type $K_r$ in a $K_r$-tiling of size at most $g(t_{I^*})$. 
Thus if we only concern about the types of $K_r$'s (that is, identify copies of $K_r$ of same type), then the number of distinct $K_r$-tilings of size at most $g(t_{I^*})$ is equal to the number of solutions of $\sum_{i=1}^{y}x_i\le g(t_{I^*})$, i.e. the number of solutions of $\sum_{i=1}^{y}x_i+x_{y+1}=g(t_{I^*})$ with $x_1, \dots, x_{y+1}\in \{0\}\cup \mathbb{N}$, which is  ${g(t_{I^*})+y \choose y}\le (g(t_{I^*})+y)^y=O(n^y)$.
For each solution $(x_1,x_2,\ldots,x_y)$ of $\sum_{i=1}^{y}x_i\le g(t_{I^*})$,
we can decide whether $G$ has a $K_r$-tiling of size $C:=\sum_{i=1}^{y}x_i$ in time $O(n^{r+1}2^{O(rC)})=O(n^{r+1}2^{O(g(t_{I^*}))})$ by Lemma \ref{lem:color coding}.
Hence the running time of this step is $O(n^y) \cdot O(n^{r+1}2^{O(g(t_{I^*}))})=n^{O(1)}2^{O(g(t_{I^*}))}$.

\noindent \textbf{Step 2.} Note that checking the slack of a $K_r$-tiling takes $O(n)$ time. Since Step 1 returns at most ${g(t_{I^*})+y \choose y}=O(n^y)$ $K_r$-tilings, the running time of this step is $O(n^y)\cdot O(n)=n^{O(1)}$. 
If there is a $K_r$-tiling $F$ obtained in Step 1 satisfying $\textbf{s}(F)|_{i}=t_{i}$ for each $i\in I^*$, then we take it as $K^1$.
Otherwise, by Lemma \ref{lem:slack}, there is no $K_r$-factor in $G$, then we halt and output the partition $(A_1^1,\ldots,A_s^1,B^1)$ of $G$ as a certificate. 

Note that the total running time of the above two steps is $n^{O(1)}2^{O(g(t_{I^*}))}$, which is $n^{O(1)}2^{O(c^{r^{r+2}})}$ for $r\ge 4$ and $n^{O(1)}2^{O(c)}$ for $r=3$.
\end{proof}


Before proving Lemma~\ref{lem:balancedness2}, we first present the following useful result.

\begin{lemma}\label{lem:matching}
If $t_i>c/0.99$ for a fixed $i\in[s]$, then $G[A_i^1]$ has a matching of size at least $xt_i$, where $x<1/(2400r\mu^{1/3})$.
\end{lemma}
\begin{proof}
It suffices to prove that there is an edge even after deleting any set $U\subseteq A_i^1$ of size $2(xt_i-1)$.
Since $|A_i^1|=\frac{n}{r}+t_i$, we have $\delta(G[A_i^1])\ge t_i-c$. On the one hand, $e(G[A_i^1])\ge \frac{1}{2}|A_i^1|\delta(G[A_i^1])\ge \frac{1}{2}(\frac{n}{r}+t_i)(t_i-c)\ge \frac{n}{2r}(t_i-c)$.
On the other hand, $2(xt_i-1)\Delta(G[A_i^1])\le 2(xt_i-1)(3\mu^{1/3}n+2c)\le 2xt_i\cdot6\mu^{1/3}n$, as $c = o(n)$. Hence there is an edge in $G[A_i^1\setminus U]$ as long as $\frac{n}{2r}(t_i-c)>2xt_i\cdot6\mu^{1/3}n$, that is,
$x<1/(2400r\mu^{1/3})$.
\end{proof}

Now we prove Lemma \ref{lem:balancedness2}.

\begin{proof}[{Proof of Lemma \ref{lem:balancedness2}.}]
Let ${A_i^1}'=A_i^1 \setminus V(K^1)$ and ${B^1}'=B^1 \setminus V(K^1)$. Let $(t_1',t_2',\ldots,t_{s+1}')$ be the slack of partition $({A_1^1}',\ldots,{A_s^1}',{B^1}')$. Then we have $t_{i}'=0$ for each $i\in I^*$  and $\sum_{i\in[s+1]\setminus I^*}t_{i}'=0$.
Note that $|K^1|\le g(t_{I^*})$.

For each $i\in [s]\setminus I^*$, we have $t_{i}'\le t_{i}+|K^1|\le t_i+g(t_{I^*})$ (a copy of $K_r$ in $K^1$ can only make the slack of $A_i^1$ increase by at most 1). 
Recall that $|t_{i}|\ge\max\{g(t_{I^*}),c/0.99\}$ by~\ref{item:I2}. 
Hence, if $t_{i}\ge \max\{g(t_{I^*}),c/0.99\}$, then $t_{i}'\le 2t_{i}$; if $t_{i}\le \min\{-g(t_{I^*}),-c/0.99\}$, then $t_{i}'\le 0$.
By Lemma \ref{lem:matching}, if $t_{i}\ge \max\{g(t_{I^*}),c/0.99\}$, then $A_{i}^1$ has a matching of size at least $(r+2)t_{i}$ and so $A_{i}^{1'}$ has a matching of size at least $2t_{i}$.
Hence if $t_{i}'>0$, then we can find a matching $M_{i}$ of size $t_{i}'$ in time $O(n^4)$ by the Edmonds algorithm.

If $s+1\notin I^*$, then we have $t_{s+1}'\le t_{s+1}+(r-s)|K^1|$. Hence, if $t_{s+1}\ge g(t_{I^*})$, then $t_{s+1}'\le (r-s+1)t_{s+1}<s(r-s+1)\mu^{2/3}n$, as $t_{s+1}<s\mu^{2/3}n$; if $t_{s+1}\le -g(t_{I^*})$, then $t_{s+1}'\le (r-s-1)g(t_{I^*})=o(n)$, as $g(t_{I^*})=o(n)$. 
By Lemma \ref{lem:Kr-s+1}, there are at least $\mu^{7/12}n^{r-s+1}$ copies of $K_{r-s+1}$ in $G[B_g^1]$, and so $G[B_g^{1'}]=G[B_g^1\setminus V(K^1)]$ contains at least $\mu^{7/12}n^{r-s+1}-r|K^1|n^{r-s}\ge (1/2)\mu^{7/12}n^{r-s+1}$ copies of $K_{r-s+1}$. 
Hence if $t_{s+1}'>0$, then we can use a greedy algorithm to find a $K_{r-s+1}$-tiling $T$ of size $t_{s+1}'$ in $G[B_g^{1'}]$ in time $O(n^{r-s+2})$, as $\mu^{2/3}<\mu^{7/12}$.

Within this proof, we may relabel the $A_i^1$'s so that $t_i'>0$ for $1\le i\le a$, $t_i'=0$ for $a+1\le i\le b$, $t_i'<0$ for $b+1\le i\le s$.
Then we shall find $\sum_{i=1}^{a}t_i'+\max\{t_{s+1}',0\}$ vertex-disjoint unbalanced copies of $K_r$ by extending $M_i$'s and $T$, which together with $K^1$ give the desired $K_r$-tiling $K$.

\begin{algorithm}
        \caption{Define $t_{i,k}'$}
        \label{alg: t_ik}
        \KwData{$t_1',t_2',\ldots,t_{s+1}'$}
        \KwResult{$t_{i,k}'$'s for each $i\in [1,a]$ and $k\in [b+1,s+1]$}
        \For{$i=1$ to $a$}{
        \For{$k=b+1$ to $s+1$}{
        \eIf{$\sum_{i'=1}^{i-1}t_{i',k}'<|t_k'|$}
         {set $t_{i,k}'=\min\{|t_k'|-\sum_{i'=1}^{i-1}t_{i',k}', t_i'-\sum_{k'=b+1}^{k-1}t_{i,k'}'\}$; \\
         }
         {set $t_{i,k}'=0$; \\
         }
        }
        }
\end{algorithm}

First, we need to decide the types of the copies of $K_r$ we need, which is done by Algorithm \ref{alg: t_ik}.
Given a sequence of input $t_1',t_2',\ldots,t_{s+1}'$ with $t_1',\dots,t_a'>0$ and $t_{b+1}',\dots,t_{s+1}'
<0$, Algorithm \ref{alg: t_ik} defines a set of numbers $t_{i,k}'$ for each $i\in [1,a]$ and $k\in [b+1,s+1]$.
If $t_{s+1}'\le 0$, then we have $t_i'=\sum_{k=b+1}^{s+1}t_{i,k}'$ for any $i\in[a]$ and $|t_k'|=\sum_{i=1}^{a}t_{i,k}'$ for any $k\in [b+1,s+1]$.
This is exactly what we want: if we extend $M_i$'s and $T$ to a $K_r$-tiling $K^2$ among which exactly $t_{i,k}'$ of them are $(i,k)$-copies (defined in Definition~\ref{def:index_vector}) for each $i\in [a]$ and $k\in [b+1,s+1]$, then $K^2$ has the same slack as $G-K^1$, that is, $G-(K^1\cup K^2)$ is balanced.
The case $t_{s+1}'>0$ can be treated similarly.

Formally, if $t_{s+1}'\le 0$, then define $t_{i,k}'$ by Algorithm \ref{alg: t_ik}. 
Note that $t_i'=\sum_{k=b+1}^{s+1}t_{i,k}'$ for any $i\in[a]$ and $|t_k'|=\sum_{i=1}^{a}t_{i,k}'$ for any $k\in [b+1,s+1]$.
If $t_{s+1}'>0$, then we instead run Algorithm \ref{alg: t_ik} with the for-loops for $i\in [a]\cup\{s+1\}$ and $k\in [b+1,s]$. In particular, when $i=s+1$, we need to check if $\sum_{i'=1}^{a}t_{i',k}'<|t_k'|$ for each $k\in [b+1,s]$. 
Note that in this case we have $t_i'=\sum_{k=b+1}^{s}t_{i,k}'$ for any $i\in[a]\cup \{s+1\}$ and $|t_k'|=\sum_{i\in [a]\cup \{s+1\}}t_{i,k}'$ for any $k\in [b+1,s]$.
The running time of Algorithm \ref{alg: t_ik} is $O(1)$.

It remains to show that there are indeed $t_{i,k}'$ vertex-disjoint $(i,k)$-copies of $K_r$ in $G-K^1$ for each $t_{i,k}'$.
Note that for each $i\in [s]$, every vertex in $A_i^1$ has at most $2\mu^{1/3}n+c+\mu^{2/3}n$ non-neighbors in $V(G)\setminus A_i^1$ by Fact \ref{fact: properties of Ai-good vtx} (1); moreover, every vertex in $B^1_g$ has at most $2\mu^{2/5}n+4\mu rn+2c$ non-neighbors in $A_i^1$ for each $i\in [s]$ by Definition \ref{def: B1-bad vtx} and Fact \ref{fact: properties of Ai-good vtx} (2), and $\delta(G[B^1_g])\ge (r-s-1)\frac{n}{r}+t_{s+1}-3s\mu^{3/5}n-c$.
Hence we can extend $M_i$ greedily to a $K_r$-tiling of size $t_i'=\sum_{k=b+1}^{s+1}t_{i,k}'$ for any $i\in[a]$, and extend $T$ greedily to a $K_r$-tiling of size $t_{s+1}'=\sum_{k=b+1}^{s}t_{s+1,k}'$ (if $t_{s+1}'>0$).
The running time of this process is $O(n^r)$. 
Let $K^2$ be the union of these $K_r$-tilings.


Let $K=K^1\cup K^2$ and $(A_1^2,\ldots,A_s^2,B^2)=(A_1^1,\ldots,A_s^1,B^1)-K$. Then $A_1^2,\ldots,A_s^2,B^2$ is a balanced partition with $|A_1^2|=\cdots=|A_s^2|=\frac{|B^2|}{r-s}=n/r-|K|\ge n/r-r^2\mu^{2/3}n$, as $t_i'<8\mu rn+4c$ for each $i\in [s]$, $t_{s+1}' <s(r-s+1)\mu^{2/3}n$ and $c=o(n)$.

The ``in particular'' part of the lemma can be easily seen from the proof. Indeed, $|B^1\setminus V(K^1)|-(r-s)(n/r-|K^1|)=t_{s+1}' > 0$, and so $K$ contains a copy of $K_r$ that has exactly $r - s + 1$ vertices of $B^1$.
\end{proof}

We finish this Section \ref{sec: Auxiliary lemmas} by proving Lemma \ref{lem:clean bad vtx-pre}.

\begin{proof}[Proof of Lemma \ref{lem:clean bad vtx-pre}.]
Recall that $d_{A_i}(v)> 2\mu^{1/3}n$ for any $i\in [s]$ and $d_{A_i}(v)< |A_i|-2\mu^{2/5}n$ for some $i \in [s]$.
Note that since $d_G(v)\ge (1-1/r)n-c$, there is at most one of $A_1,A_2,\ldots,A_s$ such that $v$ has at most $0.4n/r$ neighbors in it.
Without loss of generality, assume that $d_{A_1}(v)$ is the smallest among all $d_{A_i}(v)$.
Then $d_{A_1}(v)> 2\mu^{1/3}n$ and $d_{A_i}(v)\ge 0.4n/r$ for any $i\in [2,s]$.

We claim that there are many balanced copies of $K_r$ containing $v$ in $G$. To show this claim, we first give some basic properties:
\begin{enumerate}[label=(P\arabic*)]
  \item Since $\delta(G[B^1])\ge (r-s-1)n/r+t_{s+1}-c\ge (r-s-1)n/r-\mu ^{1/3}n$ and $|B_b^1|\le 3s\mu^{3/5}n$ (Fact~\ref{fact: properties of B1-bad vtx}), any $r-s-1$ vertices of $B^1$ have at least $\frac{n}{2r}$ common neighbors in $B^1_g$.

  \item Since $d_{A_1}(v)> 2\mu^{1/3}n$, we have $d_{A_1^1}(v)> 2\mu^{1/3}n-\mu^{2/3}n>\mu^{1/3}n$.\label{item:P5}
  
  \item For any $i\in [2,s]$, since $d_{A_i}(v)> 0.4n/r$, we have $d_{A_i^1}(v)>0.4n/r-\mu^{2/3}n$. 
  \label{item:P2}
  \item Let $u\in B^1_g$.
      By Definition \ref{def: B1-bad vtx}, $d_{A_i}(u)> |A_i|-2\mu^{2/5}n$ for any $i\in [s]$.
      Then by Fact \ref{fact: properties of Ai-good vtx} (2), $d_{A_i^1}(u)> |A_i^1|-2\mu^{2/5}n-4\mu rn-2c> |A_i^1|-3\mu^{2/5}n$ 
      for any $i\in [s]$.
  \item Given $i\in [s]$, let $w\in A_i^1$. Recall that $w$ is an $A_i$-good vertex. 
      By Definition \ref{def: Ai-good vtx}, $d_{A_i}(w)< 2\mu^{1/3}n$.
      Then for any $j\in [s]\setminus \{i\}$, we have $d_{A_j}(w)>|A_j|-2\mu^{1/3}n-c$, $d_{A_j^1}(w)>|A_j^1|-2\mu^{1/3}n-c-(4\mu rn+2c)> |A_j^1|-3\mu^{1/3}n$.
      \label{item:P4}
\end{enumerate}

Hence, we can greedily build a balanced copy of $K_r$ containing $v$, by first choosing vertices in $B_g^1$, then a vertex of $A_1^1$ and then vertices in each of $A_i^1$, $i\in [2,s]$, while each time we avoid any vertex of $W$.
Since at each time the number of common neighbors is at least $\mu^{1/3}n-(r-s-1)\cdot3\mu^{2/5}n\ge \mu^{1/3}n/2$ (indeed, achieved by the number of choices for a vertex in $A_1^1$), and $|W|\le r^3\mu^{3/5}n$, we can always choose a desired vertex and at the end of a desired copy of $K_r$ in time $O(n)$.
%
%
\end{proof}

\subsection{Proof of Lemma \ref{lem:charact}.}\label{sec: obtain K_r-s-factor}
We start with applying Lemma~\ref{lem:balancedness} to $G$ with the partition $A_1^1,\ldots,A_s^1,B^1$ and the set $I^*\subseteq [s+1]$.
If the lemma outputs the partition $A_1^1,\ldots,A_s^1,B^1$ showing that no $I^*$-admissible $K_r$-tilings exist, by Lemma~\ref{lem:slack}, $G$ has no $K_r$-factor.
So it remains to consider the case Lemma~\ref{lem:balancedness} outputs an $I^*$-admissible $K_r$-tiling $K^1$.

We first apply Lemma~\ref{lem:balancedness2} to find a $K_r$-tiling $K\supseteq K^1$ of size at most $r^2\mu^{2/3}n$ in time $O(n^r)$ such that $K^1\subseteq K$ and $V(G)\setminus V(K)$ has a partition $A_1^2,\ldots,A_s^2,B^2$ with $|A_1^2|=\cdots=|A_s^2|=\frac{|B^2|}{r-s}\ge n/r-r^2\mu^{2/3}n$, $A_i^2\subseteq A_i^1$ for each $i\in [s]$ and $B^2\subseteq B^1$. 
In particular, if $|B^1\setminus V(K^1)| - (r-s)(n/r - |K^1|) > 0$, then $K$ contains a copy of $K_r$ that has exactly $r - s + 1$ vertices of $B^1$.
Then we use Lemma~\ref{lem:clean bad vtx-pre} repeatedly to find a $K_r$-tiling $K^b$ whose members are all balanced copies of $K_r$ and contain all bad vertices of $B^1$ not covered by $K$.

Let $A_1^3,\ldots,A_s^3,B^3$ be the balanced partition of $G-(K\cup K^b)$ with $|A_1^3|=\cdots=|A_s^3|=\frac{|B^3|}{r-s}\ge n/r-r^2\mu^{2/3}n-s(2\mu^{3/5}+\mu^{2/3})n\ge n/r-3s\mu^{3/5}n$.
Note that we have
\begin{enumerate}[label=(\roman*)]
  \item for every vertex $v\in B^3\subseteq B_g^1$ and $i\in [s]$, we have $d_{A_i^3}(v)>|A_i^3|-3\mu^{2/5}n$;
  \item for any $i\in [s]$ and $u\in A_i^3$, we have $d_{A_j^3}(u)>|A_j^3|-3\mu^{1/3}n$ for any $j\in [s]\setminus \{i\}$ and $d_{B^3}(u)> |B^3|-3\mu^{1/3}n$.
\end{enumerate}
Indeed, since every vertex $v\in B_g^1$ is not bad, we have $d_{A_i}(v)>|A_i|-2\mu^{2/5}n$ and so $d_{A_i^3}(v)>|A_i^3|-2\mu^{2/5}n-(4\mu rn+2c)>|A_i^3|-3\mu^{2/5}n$ for any $i\in [s]$.
Moreover, for any $i\in [s]$ and $u\in A_i^3$, we have $d_{A_j^3}(u)>|A_j^3|-2\mu^{1/3}n-c-(4\mu rn+2c)>|A_j^3|-3\mu^{1/3}n$ for any $j\in [s]\setminus \{i\}$ and $d_{B^3}(u)>|B^3|-2\mu^{1/3}n-c-s\mu^{2/3}n\ge |B^3|-3\mu^{1/3}n$.

Properties (i) and (ii) are needed for us to apply Lemma~\ref{lem: ideal} to $G[A_1^3\cup \dots\cup A_s^3\cup B^3]$ or a slight modification of it. 

The rest of the proofs is split into several cases depending on whether $r-s=2$, $I^*=[s+1]$, $B_b^1=\emptyset$ and whether $B^1$ is disconnected with two odd components.
These simple properties can be checked easily though we need a simple Breath-First Search for the last.

\begin{claim}\label{lem:r-s neq 2}
If $r-s\neq 2$, then we can find a $K_{r}$-factor of $G$ in time $O(n^{4r^2})$.
\end{claim}

\begin{proof}

Since $(r-s)n/r-3r^2\mu^{3/5}n\le |B^3|\le (r-s)n/r$, we have $\delta(G[B^3])\ge (r-s-1)n/r-3r^2\mu^{3/5}n-c\ge \left(1-\frac{1}{r-s}-6r^3\mu^{3/5}-\frac{2rc}{n}\right)|B^3|$.
Moreover, since $B^1$ has no $\frac{\lambda}{2}$-sparse $\frac{n}{r}$-set of $G$, it follows that $B^3$ has no $\frac{\lambda}{4}$-sparse $\frac{|B^3|}{r-s}$-set.
Since $r-s\neq 2$, $\mu\ll \lambda$ and $c = o(n)$, we can apply Theorem \ref{thm:non-ext} to $G[B^3]$ and obtain a $K_{r-s}$-factor of $G[B^3]$ in time $O(n^{4r^2})$. 
Recall properties (i) and (ii).
By applying Lemma \ref{lem: ideal} with $A_i'=A_i^3$ for any $i\in [s]$ and $B'=B^3$, we can find  a $K_{r}$-factor of $G[A_1^3\cup \dots\cup A_s^3\cup B^3]$ in time $O(sn^4)$. Together with $K$ and $K^b$, this gives a $K_{r}$-factor of $G$.
\end{proof}


Now suppose $r-s=2$.
We use the Edmonds algorithm to check if $B^3$ has a perfect matching in time $O(n^4)$.
First suppose it does, and the Edmonds algorithm outputs it. Hence, by applying Lemma \ref{lem: ideal} with $A_i'=A_i^3$ for any $i\in [s]$ and $B'=B^3$, we find a $K_{r}$-factor of $G[A_1^3\cup \dots\cup A_s^3\cup B^3]$ in time $O(sn^4)$. Together with $K$ and $K^b$, we obtain a $K_{r}$-factor of $G$.

So we may assume that $B^3$ has no perfect matching.
Then by Lemma~\ref{lem:odd comp}, we conclude that $B^3$ has two odd components both of which are almost complete, denoted by $B^*$ and $B^{**}$.
We first give a simple and useful result.

\begin{prop}\label{claim3}
Let $H$ be a copy of $K_r$ in $K$. We can extend $V(H) \setminus B^1$ to a copy $H'$ of $K_r$ in time $O(n^r)$ such that $|V(H')\cap B^*|=|V(H)\cap B^1|$ or  $|V(H')\cap B^{**}|=|V(H)\cap B^1|$. 
\end{prop}
\begin{proof}
  The result holds because vertices in $V(H) \setminus B^1$ are adjacent to almost all vertices in $B^*$ (resp. $B^{**}$) and $G[B^*]$ (resp. $G[B^{**}]$) is almost complete.   
\end{proof}

\begin{claim}\label{lem: r-s=2 case1}
If $r-s = 2$ and $I^*\neq [s+1]$, then we can find a $K_r$-factor of $G$ in time $O(n^{r+1})$.
\end{claim}

\begin{proof}

Since $I^*\neq [s+1]$, we have $|A_i^1\setminus V(K^1)| - (n/r - |K^1|) > 0$ for some $i\in [s]$ or $|B^1\setminus V(K^1)| - (r-s)(n/r - |K^1|) > 0$. Suppose the former one holds. 
Then there is a large $t_i$ for some $i\in [s]\setminus I^*$, namely, $t_i  > c/0.99$.
Without loss of generality, set $i=1$, and then $t_1  > c/0.99$.
By Lemma~\ref{lem:matching} there is an edge in $A_1^3$, which can be found in $O(n^2)$. Then we extend it to a copy $H_1$ of $K_r$ with index vector $(2,1,\ldots,1,1)$ by using one vertex from $B^*$.
We next find a copy of $K_3$ in $B^{**}$, and then extend it to a copy $H_2$ of $K_r$ with index vector $(0,1,\ldots,1,3)$ and it is vertex-disjoint from $H_1$. 
The running time of the above process is $O(n^{r-2})$. 

Note that both $B^*\setminus V(H_1)$ and $B^{**}\setminus V(H_2)$ are even and so $B^3\setminus (V(H_1)\cup V(H_2))$ has a perfect matching. Moreover, 
the sizes of $A_1^3\setminus (V(H_1)\cup V(H_2)),\ldots,A_s^3\setminus (V(H_1)\cup V(H_2)),B^3\setminus (V(H_1)\cup V(H_2))$ satisfy the ratio $1:\cdots : 1:2$.
By applying Lemma \ref{lem: ideal} with $A_i'=A_i^3\setminus (V(H_1)\cup V(H_2))$ for any $i\in [s]$ and $B'=B^3\setminus (V(H_1)\cup V(H_2))$, we can find  a $K_{r}$-factor of $G[A_1'\cup \dots\cup A_s'\cup B']$ in time $O(sn^4)$. Together with $\{H_1,H_2\}$, $K$ and $K^b$, we obtain a $K_{r}$-factor of $G$.

Next, assume $|B^1\setminus V(K^1)| - (r-s)(n/r - |K^1|) > 0$.
In this case, recall that $K$ contains a copy of $K_r$, denoted by $H$, such that $|V(H)\cap B^1|=r-s+1=3$.
Now if $G[B^3\cup (V(H)\cap B^1)]$ is connected, then by Proposition~\ref{claim3}, let $H_1$ be a copy of $K_r$ that consists of $V(H)\setminus B^1$ and a triangle from $B^*$.
Then either $G[(B^3\cup V(H))\setminus V(H_1)]$ is connected or it has two even components, and in either case, it has a perfect matching.
Then we can find a $K_r$-factor of $G[A_1^3\cup \cdots \cup A_s^3\cup (B^3\cup V(H))\setminus V(H_1)]$ by Lemma~\ref{lem: ideal} similarly as above, which together with $((K\cup K^b)\setminus \{H\})\cup \{H_1\}$ yields a $K_r$-factor of $G$.
The remaining case is that $G[B^3\cup (V(H)\cap B^1)]$ is disconnected.
Without loss of generality, assume that $V(H)\cap B^1$ has neighbors only in $B^*$.
Then by Proposition~\ref{claim3}, let $H_2$ be a copy of $K_r$ that consists of $V(H)\setminus B^1$ and a triangle from $B^{**}$.
Note that $G[(B^3\cup V(H))\setminus V(H_2)]$ has two even components and we can find a $K_r$-factor of $G$ by taking the union of the $K_r$-factor returned by Lemma~\ref{lem: ideal} and $((K\cup K^b)\setminus \{H\})\cup \{H_2\}$.
\end{proof}

\begin{claim}\label{lem: r-s=2 have bad vtx}
If $r-s=2$ and $B^1_b\neq \emptyset$, then we can find a $K_r$-factor of $G$ in time $O(n^{r+1})$.
\end{claim}

\begin{proof}

Let $v\in B^1_b$. Then $v\in V(K\cup K^b)$. Recall that $|B^3|\ge 2n/r-6s\mu^{3/5}n$ and $\delta(G[B^3])\ge n/r-6s\mu^{3/5}n-c$. 
Since $v$ is a bad vertex, there exists some $i\in [s]$ such that 
$d_{A_i}(v)<|A_i|-2\mu^{2/5}n$.
Moreover, by Fact~\ref{fact: properties of Ai-good vtx}(1), $|B^3\cap A_i|\le \mu^{2/3}n$.
Then $v$ has at most $n/r+c - (2\mu^{2/5}n-\mu^{2/3}n)$ non-neighbors in $B^3$.
Therefore, for any set $Y\subseteq B^3$ of size at least $n/r-6s\mu^{3/5}n-c-2r$, we have
\begin{equation}
\label{eq:dyv}
d_Y(v)>|Y|-n/r-c+(2\mu^{2/5}n-\mu^{2/3}n)>0.
\end{equation}

Suppose $H$ is the copy of $K_r$ in $K\cup K^b$ containing $x\in [r]$ bad vertices of $B^1$, including $v$. 
We extend $V(H)\setminus B^1$ to a copy $H'$ of $K_r$ in time $O(n^r)$ such that $|V(H')\cap B^*|=|V(H)\cap B^1|$ by Proposition~\ref{claim3}.

Next we cover the $x-1$ bad vertices except $v$ by a small $K_r$-tiling. 
Indeed, we use Lemma~\ref{lem:clean bad vtx-pre} repeatedly to find a $K_r$-tiling $F$ whose members are all balanced copies of $K_r$ and contain all these bad vertices.
Let $\tilde B$ be the remaining of vertices of $B^1$, excluding $v$.
Clearly, $|\tilde B|$ is odd and $\tilde B\subseteq B_g^1$.
Note that $d_{\tilde B}(v)\ge 0.9n/r$.
If $\tilde B$ is connected, then there exists a vertex $y\in \tilde B\cap N(v)$ such that $\tilde B\setminus \{y\}$ is still connected.
Indeed, if such $y$ does not exist, then all neighbors of $v$ are cut-vertices.
In particular, $\tilde B$ has two cut-vertices, which leaves a 2-connected component $Z$ of size at most $2n/(3r)$.
Then any vertex in $Z$ has degree at most $2n/(3r)+2 < 0.99n/r$ in $\tilde B$, a contradiction.  
Then we can find a balanced copy $T$ of $K_r$ containing the edge $vy$, which is possible by~\ref{item:P5}--\ref{item:P4}. 
On the other hand, if $\tilde B$ is disconnected, then note that it has an even component $X$ and an odd component $Y$.
Recalling that $\delta(G[B^3])\ge n/r-6s\mu^{3/5}n-c$ and $|F|\le r$, we infer that $|Y|\ge n/r-6s\mu^{3/5}n-c - 2r$.
By~\eqref{eq:dyv}, we have $d_Y(v)>0$.
Then similarly we pick a balanced copy $T$ of $K_r$ containing the edge $vy$, where $y\in Y$ by ~\ref{item:P5}--\ref{item:P4}.
Therefore, in both cases, we have that $\tilde B\setminus V(T)=\tilde B\setminus \{y\}$ has a perfect matching.

By applying Lemma \ref{lem: ideal} with $A_i'=A_i^3\setminus (V(F)\cup V(T))$ for any $i\in [s]$ and $B'=\tilde B\setminus V(T)$, we can find  a $K_{r}$-factor of $G[A_1'\cup\dots\cup A_s'\cup B']$ in time $O(sn^4)$. Combining with $K\cup \{H'\}\setminus \{H\}$, $F\cup \{T\}$ and $K^b$, we obtain a $K_{r}$-factor of $G$.
\end{proof}

It remains to consider the case $r-s=2$, $I=[s+1]$ and $B^1_b=\emptyset$.
Since $I=I^*=[s+1]$, $K^1$ is indeed $[s+1]$-admissible and $K=K^1$. Moreover, since $B^1_b=\emptyset$, we have $K^b=\emptyset$. 

\begin{claim}\label{lem: r-s=2 case2}
If $r-s=2$, $I=[s+1]$, $B^1_b=\emptyset$ and $B^1$ is not disconnected with two odd components, then we can find a $K_r$-factor of $G$ in time $O(n^{r+1})$.
\end{claim}

\begin{proof}
Recall that $B^3$ is disconnected with two odd components both of which are almost complete, denoted by $B^*$ and $B^{**}$.

Let $B_1^*:=\{v\in B^1: d_{B^*}(v)\ge 0.1|B^1|\}$ and $B_1^{**}:=\{v\in B^1: d_{B^{**}}(v)\ge 0.1|B^1|\}$. Note that we can find $B_1^*$ and $B_1^{**}$ in time $O(n^{2})$.
By the minimum degree condition of $B^1$, we know $B^1 = B_1^*\cup B_1^{**}$. Moreover, $B^*\subseteq B_1^*$ and $B^{**}\subseteq B_1^{**}$.

Suppose $B^\circ:=B_1^*\cap B_1^{**}\neq \emptyset$. 
Since $B^*\cap B^{**}= \emptyset$, we have $B^\circ \subseteq B^1\setminus B^3$ and so $B^\circ \subseteq V(K)$.
Let $v$ be an arbitrary vertex in $B^\circ$ and $H$ be the copy of $K_r$ in $K$ containing $v$. We extend $V(H) \setminus B^1$ to a copy $H'$ of $K_r$ such that $|V(H')\cap B^*|=|V(H)\cap B^1|$ by Proposition \ref{claim3}. Note that $(B^3\cup V(H))\setminus V(H')$ is connected via $v$ and has a perfect matching which can be found in time $O(n^4)$. By applying Lemma \ref{lem: ideal} with $A_i'=A_i^3$ for any $i\in [s]$ and $B'=(B^3\cup V(H))\setminus V(H')$, we can find  a $K_{r}$-factor of $G[A_1'\cup \dots\cup A_s'\cup B']$ in time $O(sn^4)$. Together with $K\cup \{H'\}\setminus \{H\}$, we obtain a $K_{r}$-factor of $G$.

Next suppose $B_1^*\cap B_1^{**} = \emptyset$ and $E(G[B_1^*,B_1^{**}])\neq \emptyset$. Let $uv\in E(G[B_1^*,B_1^{**}])$ with $u\in B_1^*$ and $v\in B_1^{**}$. Since there is no edge between $B^*$ and $B^{**}$, either $u$ or $v$ is in $V(K)$. By Proposition \ref{claim3}, we swap $uv$ out. Indeed, let $H_u$ (resp. $H_v$) (if exists) be the copy of $K_r$ in $K$ containing $u$ (resp. $v$). Note that at least one of $H_u$ and $H_v$ exists. Assume that both $H_u$ and $H_v$ exist and $H_u\neq H_v$ (other cases are similar). Then we apply Proposition \ref{claim3} twice and obtain $H_u',H_v'$. Note that $(B^3\cup V(H_u)\cup V(H_v))\setminus (V(H_u')\cup V(H_v'))$ is connected via $uv$ and has a perfect matching  which can be found in time $O(n^4)$.
By applying Lemma \ref{lem: ideal} with $A_i'=A_i^3$ for any $i\in [s]$ and $B'=(B^3\cup V(H_u)\cup V(H_v))\setminus (V(H_u')\cup V(H_v'))$, we can find  a $K_{r}$-factor of $G[A_1'\cup \dots\cup A_s'\cup B']$ in time $O(sn^4)$. Combining with $K\cup \{H_u', H_v'\}\setminus \{H_u,H_v\}$, we obtain a $K_{r}$-factor of $G$.

Finally, suppose $B^1$ is disconnected with two components $B_1^*$ and $B_1^{**}$. By the assumption of the claim, at least one of $|B_1^*|$ and $|B_1^{**}|$ is even. Without loss of generality, assume that $|B_1^*|$ is even. Since $|B^*|$ is odd, there is a copy $H$ of $K_r$ in $K$ such that $|V(H)\cap B_1^*|$ is odd. By Proposition~\ref{claim3}, we can extend $V(H) \setminus B_1^*$ to a copy $H'$ of $K_r$ such that $|V(H')\cap B^{**}|=|V(H)\cap B_1^*|$. Now note that both $|(B^*\cup V(H))\setminus V(H')|$ and $|B^{**}\setminus V(H')|$ are even, and so $(B^3\cup V(H))\setminus V(H')$ induces a perfect matching  which can be found in time $O(n^4)$. By applying Lemma \ref{lem: ideal} with $A_i'=A_i^3$ for any $i\in [s]$ and $B'=(B^3\cup V(H))\setminus V(H')$, we find a $K_{r}$-factor of $G[A_1'\cup \dots\cup A_s'\cup B']$ in time $O(sn^4)$. Together with $(K\cup \{H'\})\setminus \{H\}$, we obtain a $K_{r}$-factor of $G$.
\end{proof}

Therefore, the proof of (1) in the lemma is completed, and it remains to consider the case that $B^1$ is disconnected with two odd components.
We need to show that
$G$ has a $K_r$-factor if and only if $G$ has a strong admissible $K_r$-tiling. Moreover, the existence of such $K_r$-tiling in $G$ can be checked by an algorithm that runs in time $n^{O(1)}2^{O(c^{{r}^{r+2}})}$ for $r\ge 4$ (resp. $n^{O(1)}2^{O(c)}$ for $r=3$), and if $G$ has a $K_r$-factor, it also outputs a $K_r$-factor.

Let $X,Y$ be the two odd components of $B^1$. By the minimum degree condition of $G[B^1]$, we know that $X,Y$ are almost complete and have order roughly $n/r$. Note that if we remove a strong admissible $F$ from $G$, then $|X\setminus V(F)|-|Y\setminus V(F)|$ is even because $|B^1\setminus V(F)|$ would be even.

To check whether $G$ has a strong admissible $K_r$-tiling, we fix a copy $H$ of $K_r$ with an odd number of vertices of $B^1$ and use Lemma~\ref{lem:color coding} to search for a $K_r$-tiling $F$ in $G\setminus V(H)$ such that $F\cup \{H\}$ is $[s+1]$-admissible and $|F\cup \{H\}|\le g(t_{[s+1]})+2$. 
By similar arguments in the proof of Lemma \ref{lem:balancedness}, the color-coding search takes time $n^{O(1)}2^{O(g(t_{[s+1]}))}$.
Since there are at most $n^r$ such $H$, the running time of this process is $O(n^r)\cdot n^{O(1)}2^{O(g(t_{[s+1]}))}= n^{O(1)}2^{O(g(t_{[s+1]}))}$,
which is $n^{O(1)}2^{O(c^{r^{r+2}})}$ for $r\ge 4$ and $n^{O(1)}2^{O(c)}$ for $r=3$.

Suppose we find a strong admissible $K_r$-tiling $F$. Let $H$ be the copy of $K_r$ in $F$ with an odd number of vertices of $B^1$. 
Without loss of generality, we assume that $V(H)\cap B^1\subseteq X$. We extend $V(H)\setminus B^1$ to a copy $H'$ of $K_r$ such that $V(H')\cap B^1\subseteq Y$. Let $F'=(F\cup \{H'\}) \setminus \{H\}$. 
Note that exactly one of $B^1\setminus V(F)$ and $B^1\setminus V(F')$ has two even components. 
Without loss of generality, we assume that $B^1\setminus V(F')$ has two even components and so it can induce a perfect matching. By applying Lemma \ref{lem: ideal} with $A_i'=A_i^1\setminus V(F')$ for any $i\in [s]$ and $B'=B^1\setminus V(F')$, we can find a $K_{r}$-factor of $G[A_1'\cup \dots\cup A_s'\cup B']$ in time $O(sn^4)$. Together with $F'$, we obtain a $K_{r}$-factor of $G$.


Now suppose $G$ has no strong admissible $K_r$-tiling and we need to show that $G$ has no $K_r$-factor. 
We may assume that $G$ has an $[s+1]$-admissible $K_r$-tiling $F$ of size at most $g(t_{[s+1]})$ -- otherwise, by Lemma \ref{lem:slack}, $G$ has no $K_r$-factor. 
Since $F$ is not a strong admissible $K_r$-tiling, every copy of $K_r$ in $F$ has an even number of vertices of $B^1$. 
Moreover, we claim that there is no edge in $A_i^1\setminus V(F)$ for each $i\in [s]$. Suppose not and without loss of generality, $A_1^1\setminus V(F)$ has an edge. Then we can find two disjoint copies of $K_r$ with index vector $(2,1,\ldots,1,1)$ and $(0,1,\ldots,1,3)$, respectively. 
The union of $F$ and these two copies of $K_r$ is a strong admissible $K_r$-tiling, a contradiction.   
This implies that $G-F$ contains $s$ independent sets $A_i^1\setminus V(F)$ for $i\in [s]$ and $B^1\setminus V(F)$, which consists of two odd components.
Then note that $G-F$ does not induce a $K_r$-factor.
Indeed, no copy of $K_r$ in $G-F$ contains more than $r-2$ vertices from $\bigcup_{i\in [r-2]}A_i^1\setminus V(F)$, and since it is balanced, a $K_r$-factor of $G-F$ (if it exists) must only contain balanced copies of $K_r$.
This would induce a perfect matching in $B^1\setminus V(F)$, which does not exist.
Since this holds for any such $F$ of size at most $g(t_{[s+1]})$, by Lemma~\ref{lem:slack}, $G$ has no $K_r$-factor.

This completes the proof of part (2) and also the proof of Lemma~\ref{lem:charact}.

\subsection{When $s=r$}\label{sec: s=r}
At last, we study the case $s=r$.
In this case, $B$ is an empty set. 
We define the partition $(A_1^1,\ldots,A_r^1,B^1)$ in the same way as in other cases and note that by Fact~\ref{fact: properties of Ai-good vtx}(1) $|B^1|\le r\mu^{2/3}n$. 
Then we assign vertices in $B^1$ to those $A_j^1$'s with negative slacks such that the sizes of these new $A_j^1$'s are all at most $n/r$. 
Let $(A_1^2,\ldots,A_r^2)$ be the resulting partition and note that the vertices of $B^1$ may be distributed among all parts.
Let $t_i:=|A_i^2| - n/r$ for $i\in [r]$.
As in the previous case, we first sort $t_1,t_2,\ldots,t_{r},\frac{c}{0.99},-\frac{c}{0.99}$ by value and obtain that $t_{\tau(1)}\ge\cdots\ge t_{\tau(i_1)}>\frac{c}{0.99}\ge t_{\tau(i_1+1)}\ge \cdots \ge t_{\tau(i_2)}\ge -\frac{c}{0.99} > t_{\tau(i_2+1)}\ge \cdots\ge t_{\tau(r)}$ in time $O(r^2)$.
Then we use Algorithm \ref{alg:t_i} with $r'=r$ on $t_{\tau(1)},\ldots,t_{\tau(r)}$ and output $I=[i_1',i_2']$ in time $O(r)$, with properties~\ref{item:I1}--\ref{item:I4}.
Let $I^*:=\{\tau(i): i\in I\}$.

The characterization result in this case is the following.
\begin{lemma}
[Characterization for $s=r$]
\label{lem:charact2}
Given an $n$-vertex graph $G$ with $\delta(G)\ge (1-1/r)n-c$, let $(A_1^1,\dots,A_s^1,B^1)$ be the partition defined above in this section and suppose $s=r$.
Let $(A_1^2,\ldots,A_r^2)$ and $I^*\subseteq [r]$ be defined as above.
Then $G$ has a $K_r$-factor if and only if $G$ contains an $I^*$-admissible $K_r$-tiling.
Moreover, the existence of such $K_r$-tiling in $G$ can be checked by an algorithm that runs in time $n^{O(1)}2^{O(c^{{r}^{r+2}})}$ for $r\ge 4$ (resp. $n^{O(1)}2^{O(c)}$ for $r=3$), and if $G$ has a $K_r$-factor, it also outputs a $K_r$-factor.
\end{lemma}

\begin{proof}
We start with applying Lemma~\ref{lem:balancedness} to search for an $I^*$-admissible $K_r$-tiling $K_1$ in time $n^{O(1)}2^{O(c^{{r}^{r+2}})}$ for $r\ge 4$ (resp. $n^{O(1)}2^{O(c)}$ for $r=3$).
If it does not output an $I^*$-admissible $K_r$-tiling, then by Lemma~\ref{lem:slack}, $G$ has no $K_r$-factor (and the partition $(A_1^2,\ldots,A_r^2)$ is a certificate).
So we may assume that it outputs an $I^*$-admissible $K_r$-tiling $K_1$.
Then we apply Lemma~\ref{lem:balancedness2} to find an $[r]$-admissible $K_r$-tiling $K\supseteq K^1$ of size at most $r^2\mu^{2/3}n$ such that $K^1\subseteq K$ and $V(G)\setminus V(K)$ has a partition $A_1^3,\ldots,A_r^3$ with $|A_1^3|=\cdots=|A_r^3|\ge n/r-r^2\mu^{2/3}n$, $A_i^3\subseteq A_i^2$ for each $i\in [r]$.

Next we find a $K_r$-tiling $K^b$ to cover all vertices of $B^1$ not covered by $K$.
Let $v\in B^1\cap A_{j}^3$ for some $j\in [r]$ and note that $d_{A_i}(v)\ge 2\mu^{1/3}n$ for every $i\in [r]$.
Without loss of generality, assume that $d_{A_1}(v)$ is the smallest among all $d_{A_i}(v)$.
By the minimum degree condition, we have $d_{A_1}(v)\ge 2\mu^{1/3}n$ and $d_{A_i}(v)\ge 0.4n/r$ for $i\in [2,r]$.
By~\ref{item:P2} and~\ref{item:P4}, we can greedily find a copy of $K_r$ containing $v$ and one vertex from each of $A_i^3$, $i\in [r]\setminus\{j\}$.
Moreover, as $|B^1|\le r\mu^{2/3}n$ and $|K|\le r^2\mu^{2/3}n$, one can greedily choose a $K_r$-tiling of size $|B^1\setminus V(K)|$ that covers all vertices of $B^1\setminus V(K)$.
Denote this $K_r$-tiling by $K^b$ and the remaining vertex partition by $(A_1^4,\ldots,A_r^4)$.

Finally, we apply Lemma \ref{lem: ideal} with $A_i'=A_i^4$ for any $i\in [r]$ and $B'=\emptyset$, and find  a $K_{r}$-factor of $G[A_1^4,\dots,A_r^4]$ in time $O(rn^4)$. 
Together with $K$ and $K^b$, we obtain a $K_{r}$-factor of $G$. 
\end{proof}

\section{Conclusion}

Lemmas~\ref{lem:charact} and~\ref{lem:charact2} together detail the hidden information of Lemma~\ref{lem:char} when $0<s\le r$.
Indeed, one can set the family $\mathcal F$ in Lemma~\ref{lem:char} to be either all $I^*$-admissible $K_r$-tilings or all strong admissible $K_r$-tilings, according to the partition $(A_1^1,\dots,A_s^1,B^1)$ and the sets $I$, $B_b^1$ and $B^1$.
On the other hand, by Theorem~\ref{thm:non-ext}, in the non-extremal case (that is, when Algorithm~\ref{alg:SparseSet} returns $s=0$), one can set $\mathcal F = \{\emptyset\}$.
Therefore, Lemmas~\ref{lem:charact},~\ref{lem:charact2} and Theorem~\ref{thm:non-ext} together imply Lemma~\ref{lem:char} (and thus also Theorem~\ref{thm:main}).

\bibliographystyle{abbrv}
\bibliography{Bibref}


\end{document}